\documentclass[12pt]{article}
\usepackage{amssymb,amsthm,amsmath,amscd}

\newtheorem{thm}{Theorem}[section]
\newtheorem{lem}[thm]{Lemma}
\newtheorem{prop}[thm]{Proposition}

\newtheorem*{theo}{Theorem}

\theoremstyle{definition}
\newtheorem{df}[thm]{Definition}
\newtheorem{rem}[thm]{Remark}
\newtheorem{ex}[thm]{Example}

\numberwithin{equation}{section}

\newcommand{\Homeo}{\operatorname{Homeo}}
\newcommand{\graph}{\operatorname{graph}}

\renewcommand{\phi}{\varphi}

\newcommand{\N}{\mathbb{N}}
\newcommand{\Z}{\mathbb{Z}}

\title{The absorption theorem for \\ affable equivalence relations}
\author{Thierry Giordano 
\thanks{Supported in part by a grant from NSERC, Canada} \\
Department of Mathematics and Statistics \\
University of Ottawa \\
585 King Edward, Ottawa, Ontario, Canada K1N 6N5 
\and
Hiroki Matui 
\thanks{Supported in part by a grant 
from the Japan Society for the Promotion of Science} \\
Graduate School of Science \\
Chiba University \\
1-33 Yayoi-cho, Inage-ku, Chiba 263-8522, Japan 
\and
Ian F. Putnam 
\thanks{Supported in part by a grant from NSERC, Canada} \\
Department of Mathematics and Statistics \\
University of Victoria \\
Victoria, B.C., Canada V8W 3P4 
\and 
Christian F. Skau 
\thanks{Supported in part by the Norwegian Research Council} \\
Department of Mathematical Sciences \\
Norwegian University of Science and Technology (NTNU) \\
N-7034 Trondheim, Norway}

\date{}

\begin{document}
\maketitle

\begin{abstract}
We prove a result about 
extension of a minimal AF-equivalence relation $R$ on the Cantor set $X$, 
the extension being `small' in the sense that 
we modify $R$ on a thin closed subset $Y$ of $X$. 
We show that 
the resulting extended equivalence relation $S$ is orbit equivalent 
to the original $R$, and so, in particular, $S$ is affable. 
Even in the simplest case---when $Y$ is a finite set---this result is 
highly non-trivial. 
The result itself---called the absorption theorem---is 
a powerful and crucial tool 
for the study of the orbit structure of minimal $\Z^n$-actions 
on the Cantor set, see Remark \ref{addition}. 
The absorption theorem is 
a significant generalization of the main theorem proved in \cite{GPS2}. 
However, we shall need a few key results from \cite{GPS2} 
in order to prove the absorption theorem. 
\end{abstract}

\section{Introduction}

We introduce some basic definitions 
as well as relevant notation and terminology, 
and we refer to \cite{GPS2} as a general reference 
for background and more details. 
Throughout this paper 
we will let $X$, $Y$ or $Z$ denote 
compact, metrizable and zero-dimensional spaces, i.e. compact spaces 
which have countable bases consisting of closed-open (clopen) subsets. 
Equivalently, 
the spaces are compact, metrizable and totally disconnected spaces. 
In particular, if a space does not have isolated points, 
it is homeomorphic to the (unique) Cantor set. 
We will study equivalence relations, 
denoted by $R$, $S$, $K$, on these spaces that are countable, 
i.e. all the equivalence classes are countable (including finite). 

Let $R\subset X\times X$ be a countable equivalence relation on $X$, 
and let $[x]_R$ denote the (countable) $R$-equivalence class, 
$\{y\in X\mid(x,y)\in R\}$, of $x\in X$. 
We say that $R$ is \emph{minimal}, 
if all the $R$-equivalence classes are dense in $X$. 
$R$ has a natural \emph{groupoid structure}. 
Specifically, if $(x,y),(y,z)\in R$, then 
the product of this \emph{composable} pair is defined by 
\[ (x,y)\cdot(y,z)=(x,z). \]
The inverse of $(x,y)\in R$ is defined to be $(x,y)^{-1}=(y,x)$. 
Let $R$ be given a Hausdorff, locally compact and second countable 
(equivalently, metrizable) topology $\mathcal{T}$, 
so that the product of composable pairs (with the relative topology 
from the product topology on $R\times R$) is continuous. 
Also, the inverse map is required to be a homeomorphism on $R$. 
We say that 
$(R,\mathcal{T})$ is a locally compact (principal) groupoid. 
The range map $r:R\to X$ is defined by $r((x,y))=x$, 
and the source map $s:R\to X$ is defined by $s((x,y))=y$, 
both maps being surjective. 

\begin{df}[\'Etale equivalence relation]
The locally compact groupoid $(R,\mathcal{T})$ is \emph{\'etale}, 
if $r:R\to X$ is a local homeomorphism, i.e. 
for every $(x,y)\in R$ there exists an open neighbourhood 
$U^{(x,y)}\in\mathcal{T}$ of $(x,y)$ such that 
$r(U^{(x,y)})$ is open in $X$ and 
$r:U^{(x,y)}\to r(U^{(x,y)})$ is a homeomorphism. 
\end{df}

\begin{rem}
Clearly $r$ is an open map, 
and one may choose $U^{(x,y)}$ to be a clopen set 
(and so $r(U^{(x,y)})$ is a clopen subset of $X$). 
One thus gets that $(R,\mathcal{T})$ is 
a locally compact, metrizable and zero-dimensional space. 
Also, $r$ being a local homeomorphism implies that 
$s$ is a local homeomorphism as well. 
Occasionally we will refer to the local homeomorphism condition 
as the \emph{\'etale condition}, 
and to $U^{(x,y)}$ as an \emph{\'etale neighbourhood} (around $(x,y)$). 
It is noteworthy that 
only rarely will the topology $\mathcal{T}$ on $R\subset X\times X$ 
be the relative topology $\mathcal{T}_{\text{rel}}$ from $X\times X$. 
In general, 
$\mathcal{T}$ is a finer topology than $\mathcal{T}_{\text{rel}}$. 
For convenience, we will sometimes write $R$ for $(R,\mathcal{T})$ 
when the topology $\mathcal{T}$ is understood from the context. 
\end{rem}

It is a fact that if $R$ is \'etale, then 
the diagonal of $R$, $\Delta=\Delta_X(=\{(x,x)\mid x\in X\})$, 
is homeomorphic to $X$ via the map $(x,x)\mapsto x$. 
We will often make the identification between $\Delta$ and $X$. 
Furthermore, $\Delta$ is an open subset of $R$. 
Also, $R$ admits an (essentially) unique left Haar system 
consisting of counting measures. 
(See \cite{P} for this. 
We shall not need this last fact in this paper.) 

\begin{df}[Isomorphism and orbit equivalence]\label{OE}
Let $(R_1,\mathcal{T}_1)$ and $(R_2,\mathcal{T}_2)$ be 
two \'etale equivalence relations on $X_1$ and $X_2$, respectively. 
$R_1$ is \emph{isomorphic} to $R_2$---we will write $R_1\cong R_2$---if 
there exists a homeomorphism $F:X_1\to X_2$ satisfying the following: 
\begin{enumerate}
\item $(x,y)\in R_1\Leftrightarrow(F(x),F(y))\in R_2$. 
\item $F\times F:R_1\to R_2$ is a homeomorphism, 
where $F\times F((x,y))=(F(x),F(y))$ for $(x,y)\in R_1$. 
\end{enumerate}
We say that $F$ \emph{implements} an isomorphism between $R_1$ and $R_2$. 

We say that $R_1$ is \emph{orbit equivalent} to $R_2$ 
if (i) is satisfied, and we call $F$ an \emph{orbit map} in this case. 
(The term \emph{orbit equivalence} is motivated by the important example of 
\'etale equivalence relations coming from group actions, 
where equivalence classes coincide with orbits (see below).) 
\end{df}

There is a notion of invariant probability measure associated to 
an \'etale equivalence relation $(R,\mathcal{T})$ on $X$. 
In fact, if $(x,y)\in R$, 
there exists a clopen neighbourhood $U^{(x,y)}\in\mathcal{T}$ of $(x,y)$ 
such that both $r:U^{(x,y)}\to r(U^{(x,y)})=A$ and 
$s:U^{(x,y)}\to s(U^{(x,y)})=B$ are homeomorphism, 
with $A$ a clopen neighbourhood of $x\in X$, 
and $B$ a clopen neighbourhood of $y\in X$. 
The map $\gamma=s\circ r^{-1}:A\to B$ is a homeomorphism such that 
$\graph(\gamma)=\{(x,\gamma(x))\mid x\in A\}\subset R$. 
The triple $(A,\gamma,B)$ is called a (local) \emph{graph} in $R$, 
and by obvious identifications (in fact, $U^{(x,y)}=\graph(\gamma)$) 
the family of such graphs form a basis for $(R,\mathcal{T})$. 
Let $\mu$ be a probability measure on $X$. 
We say that $\mu$ is \emph{$R$-invariant}, 
if $\mu(A)=\mu(B)$ for every graph $(A,\gamma,B)$ in $R$. 
If $(R_G,\mathcal{T}_G)$ is the \'etale equivalence relation 
associated with the free action of the countable group $G$ 
acting as homeomorphisms on $X$ (see Example \ref{grp} below), then 
$\mu$ is $R_G$-invariant if and only if $\mu$ is \emph{$G$-invariant}, 
i.e. $\mu(A)=\mu(g(A))$ for all Borel sets $A\subset X$, and all $g\in G$. 
Note that if $G$ is an amenable group 
there exist $G$-invariant, and hence $R_G$-invariant, probability measures. 
We remark that 
if $(R_1,\mathcal{T}_1)$ is orbit equivalent to $(R_2,\mathcal{T}_2)$ 
via the orbit map $F:X_1\to X_2$, then 
$F$ maps the set of $R_1$-invariant probability measures 
bijectively onto the set of $R_2$-invariant probability measures. 

\begin{ex}\label{grp}
Let $G$ be a countable discrete group acting freely 
(i.e. $gx=x$ for some $x\in X$, $g\in G$, implies $g=e$ 
(the identity of the group)) as homeomorphisms on $X$. 
Let 
\[ R_G=\{(x,gx)\mid x\in X,g\in G\}\subset X\times X, \]
i.e. the $R_G$-equivalence classes are simply the $G$-orbits. 
We topologize $R_G$ 
by transferring the product topology on $X\times G$ to $R_G$ 
via the map $(x,g)\mapsto (x,gx)$, 
which is a bijection since $G$ acts freely. 
With this topology $\mathcal{T}_G$ 
we get that $(R_G,\mathcal{T}_G)$ is an \'etale equivalence relation. 
Observe that if $G$ is a finite group, then $R_G$ is compact. 
\end{ex}

\section{AF and AF-able (affable) equivalence relations}

Let CEER be the acronym for \emph{compact \'etale equivalence relation}. 
We have the following general result about CEERs, 
cf. \cite[Proposition 3.2]{GPS2}. 

\begin{prop}\label{ceer}
Let $(R,\mathcal{T})$ be a CEER on $X$, 
where $X$ is a compact, metrizable and zero-dimensional space. 
Let $X\times X$ be given the product topology. 
\begin{enumerate}
\item $\mathcal{T}$ is the relative topology from $X\times X$. 
\item $R$ is a closed subset of $X\times X$ and 
the quotient topology of the quotient space $X/R$ is Hausdorff. 
\item $R$ is uniformly finite, that is, 
there is a natural number $N$ such that 
the number $\#([x]_R)$ of elements in $[x]_R$ is 
less than or equal to $N$ for all $x\in X$. 
\end{enumerate}
\end{prop}

In \cite{GPS2} the structure of a CEER, 
$(R,\mathcal{T})$ on $X$, is described. 
Figure 1 illustrates how the structure looks like: 
$X$ is decomposed into a finite number of $m$ disjoint clopen towers 
$T_1,T_2,\dots,T_m$, 
each of these consisting of finitely many disjoint clopen sets. 
The equivalence classes of $R$ are represented in Figure 1 
as the family of sets consisting of points 
lying on the same vertical line in each tower. 
(In the figure we have marked the equivalence class $[x_1]_R$ 
of a point $x_1\in T_1$. 
We also show the graph picture 
associated to the tower $T_m$ of height three.) 

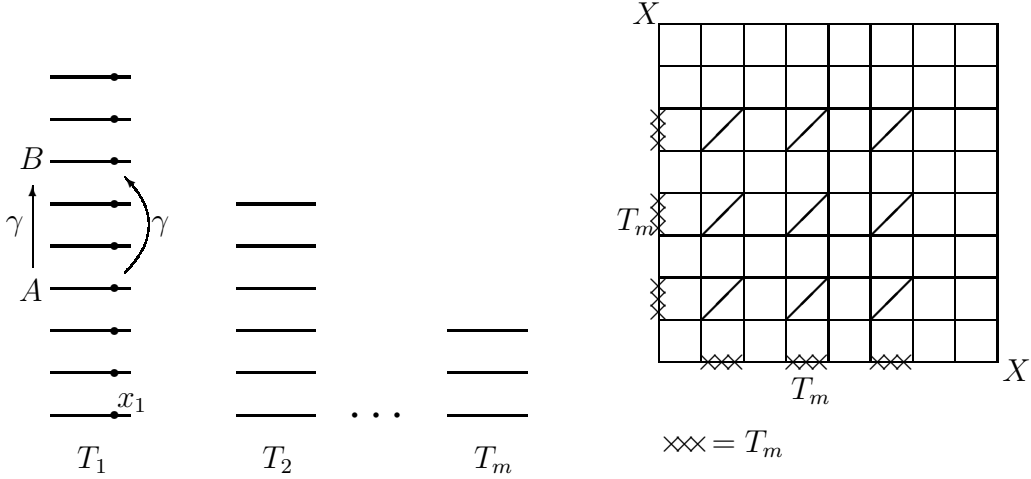
\begin{figure}
\begin{center}
\begin{picture}(400,200)(0,20)

\thicklines\multiput(20,50)(0,16){9}{\line(1,0){30}}\thinlines
\put(30,30){$T_1$}
\put(8,93){$A$}
\put(8,143){$B$}
\put(13,106){\vector(0,1){31}}
\put(3,120){$\gamma$}
\multiput(44,50)(0,16){9}{\circle*{3}}
\put(45,53){$x_1$}
\qbezier(48,104)(66,122)(48,140)
\put(50,138){\vector(-1,1){2}}
\put(58,120){$\gamma$}

\thicklines\multiput(90,50)(0,16){6}{\line(1,0){30}}\thinlines
\put(100,30){$T_2$}

\multiput(135,50)(8,0){3}{\circle*{2}}

\thicklines\multiput(170,50)(0,16){3}{\line(1,0){30}}\thinlines
\put(180,30){$T_m$}

\multiput(250,70)(0,16){9}{\line(1,0){128}}
\multiput(250,70)(16,0){9}{\line(0,1){128}}
\thicklines\multiput(266,86)(0,32){3}{\line(1,1){16}}
\multiput(298,86)(0,32){3}{\line(1,1){16}}
\multiput(330,86)(0,32){3}{\line(1,1){16}}\thinlines

\put(264,67){$\times$}\put(269,67){$\times$}\put(274,67){$\times$}
\put(296,67){$\times$}\put(301,67){$\times$}\put(306,67){$\times$}
\put(328,67){$\times$}\put(333,67){$\times$}\put(338,67){$\times$}
\put(245,86){$\times$}\put(245,91){$\times$}\put(245,96){$\times$}
\put(245,118){$\times$}\put(245,123){$\times$}\put(245,128){$\times$}
\put(245,150){$\times$}\put(245,155){$\times$}\put(245,160){$\times$}

\put(380,63){$X$}\put(240,198){$X$}
\put(300,56){$T_m$}\put(234,120){$T_m$}

\put(250,36){$\times$}\put(255,36){$\times$}\put(260,36){$\times$}
\put(270,36){$=T_m$}

\end{picture}
\caption{Illustration of the groupoid partition of 
a CEER; $R\subset X\times X$}
\end{center}
\end{figure}

Figure 1 also illustrates a concept 
that will play an important role in the sequel, namely 
a very special clopen partition of $R$, 
which we will refer to as a \emph{groupoid partition}. 
Let $A$ and $B$ be two (clopen) floors in the same tower, 
say the tower $T_1$. 
There is a homeomorphism $\gamma:A\to B$ such that 
$\graph(\gamma)=\{(x,\gamma(x))\mid x\in A\}\subset R$. 
Let $\mathcal{O}'$ be the (finite) clopen partition of $R$ 
consisting of the set of these graphs $\gamma$, 
and let $\mathcal{O}$ be the associated (finite) clopen partition of $X$ 
(which we identify with the diagonal $\Delta=\Delta_X$). 
So $\mathcal{O}=\mathcal{O}'\cap\Delta$, 
which means that $A\in\mathcal{O}$ 
if $A=B$ and $\gamma:A\to B$ is the identity map. 
The properties of $\mathcal{O}'$ are as follows 
(where we define $U\cdot V$ for subsets $U,V$ of $R$ to be 
$U\cdot V=\{(x,z)\mid (x,y)\in U, (y,z)\in V
\text{ for some }y\in X\}$): 
\begin{enumerate}
\item $\mathcal{O}'$ is a finite clopen partition of $R$ 
finer than $\{\Delta,R\setminus\Delta\}$. 
\item For all $U\in\mathcal{O}'$, 
the maps $r,s:U\to X$ are homeomorphisms onto their respective images, 
and if $U\subset R\setminus\Delta$, then $r(U)\cap s(U)=\emptyset$. 
\item For all $U,V\in\mathcal{O}'$, 
we have $U\cdot V=\emptyset$ or $U\cdot V\in\mathcal{O}'$. 
Also, $U^{-1}(=\{(y,x)\mid(x,y)\in U\})$ is in $\mathcal{O}'$ 
for every $U$ in $\mathcal{O}'$. 
\item With $\mathcal{O}'^{(2)}=\{(U,V)\mid U,V\in\mathcal{O}', \ 
U\cdot V\neq\emptyset\}$, 
define $(U,V)\in\mathcal{O}'^{(2)}\mapsto U\cdot V\in\mathcal{O}'$. 
\end{enumerate}
Then $\mathcal{O}'$ has a principal groupoid structure 
with unit space equal to $\{U\in\mathcal{O}'\mid U\subset\Delta\}$ 
(which clearly may be identified with $\mathcal{O}$). 
Hence the name \emph{groupoid partition} for $\mathcal{O}'$. 
(Note that if we think of $U$ and $V$ as maps, then 
$U\cdot V$ means first applying the map $U$ and then the map $V$.) 

Note that if we define the equivalence relation 
$\sim_{\mathcal{O}'}$ on $\mathcal{O}$ by $A\sim_{\mathcal{O}'}B$ 
if there exists $U\in\mathcal{O}'$ such that 
$U^{-1}\cdot U=A$, $U\cdot U^{-1}=B$, 
then the equivalence classes $[\cdot]_{\mathcal{O}'}$ are exactly 
the towers in Figure 1. 
The heights of the various towers $T_1,T_2,\dots,T_m$ in Figure 1 
are not necessarily distinct. 
All the groupoid partitions of $R$ finer than the one shown in Figure 1 
are obtained by vertically dividing the various towers $T_1,T_2,\dots,T_m$ 
(by clopen sets) in an obvious way. 

The proof of the following proposition can be found 
in \cite[Lemma 3.4, Corollary 3.5]{GPS2}. 

\begin{prop}\label{groupoid}
Let $(R,\mathcal{T})$ be a CEER on $X$, and let $\mathcal{V}'$ and 
$\mathcal{V}$ be (finite) clopen partitions of $R$ and $X$, respectively. 
There exists a groupoid partition $\mathcal{O}'$ of $R$ 
which is finer than $\mathcal{V}'$, and such that 
$\mathcal{O}=\mathcal{O}'\cap\Delta$ is a clopen partition of $X$ 
that is finer than $\mathcal{V}$. 
\end{prop}

\begin{df}[AF and AF-able (affable) equivalence relations]
Let $\{(R_n,\mathcal{T}_n)\}_{n=0}^\infty$ be an ascending sequence of 
CEERs on $X$ (compact, metrizable, zero-dimensional), that is, 
$R_n\subset R_{n+1}$ and $R_n\in\mathcal{T}_{n+1}$ 
(i.e. $R_n$ is open in $R_{n+1}$) for $n=0,1,2,\dots$, 
where we set $R_0=\Delta_X(\cong X)$, 
$\mathcal{T}_0$ being the topology on $X$. 
Let $(R,\mathcal{T})$ be the \emph{inductive limit} of 
$\{(R_n,\mathcal{T}_n)\}$ with the inductive limit topology $\mathcal{T}$, 
i.e. $R=\bigcup_{n=0}^\infty R_n$ and 
$U\in\mathcal{T}$ if $U\cap R_n\in\mathcal{T}_n$ for any $n$. 
In particular, $R_n$ is an open subset of $R$ for all $n$. 
We say that $(R,\mathcal{T})$ is an \emph{AF-equivalence relation} on $X$, 
and we use the notation 
$\displaystyle(R,\mathcal{T})=\lim_{\longrightarrow}(R_n,\mathcal{T}_n)$. 
We say that 
an equivalence relation $S$ on $X$ is \emph{AF-able (affable)} 
if it can be given a topology making it an AF-equivalence relation. 
(Note that this is the same as to say that 
$S$ is orbit equivalent to an AF-equivalence relation, 
cf. Definition \ref{OE}.)
\end{df}

\begin{rem}
One can prove that $(R,\mathcal{T})$ is an AF-equivalence relation 
if and only if $(R,\mathcal{T})$ is the inductive limit of 
an ascending sequence $\{(R_n,\mathcal{T}_n)\}_{n=0}^\infty$, 
where all the $(R_n,\mathcal{T}_n)$ are \'etale and \emph{finite} 
(i.e. the $R_n$-equivalence classes are finite) equivalence relations, 
not necessarily CEERs, cf. \cite{M}. 
This fact highlights the analogy 
between AF-equivalence relations in the topological setting 
with the so-called \emph{hyperfinite} equivalence relations 
in the Borel and measure-theoretic setting. 

It can be shown that 
the condition that $R_n$ is open in $R_{n+1}$ is superfluous 
when $R_n$ and $R_{n+1}$ are CEERs 
(see the comment right after Definition 3.7 of \cite{GPS2}). 
\end{rem}

We will assume some familiarity 
with the notion of a \emph{Bratteli diagram} 
(cf. \cite{GPS2} for details). 
We remind the reader of the notation we will use. 
Let $(V,E)$ be a Bratteli diagram, 
where $V$ is the vertex set and $E$ is the edge set, and 
where $V$, respectively $E$, can be written 
as a countable disjoint union of finite non-empty sets: 
\[ V=V_0\cup V_1\cup V_2\cup\dots
\quad\text{and}\quad E=E_1\cup E_2\cup\dots \]
with the following property: 
an edge $e$ in $E_n$ connects 
a vertex $v$ in $V_{n-1}$ to a vertex $w$ in $V_n$. 
We write $i(e)=v$ and $t(e)=w$, 
where we call $i$ the source (or initial) map and 
$t$ the range (or terminal) map. 
So a Bratteli diagram has a natural grading, 
and we will say that $V_n$ is the vertex set at level $n$. 
We require that $i^{-1}(v)\neq\emptyset$ for all $v\in V$ 
and $t^{-1}(v)\neq\emptyset$ for all $v\in V\setminus V_0$. 
We also want our Bratteli diagram to be \emph{standard}, 
i.e. $V_0=\{v_0\}$ is a one-point set. 
In the sequel all our Bratteli diagrams are assumed to be standard, 
so we drop the term `standard'. 
Let 
\[ X_{(V,E)}=\{(e_1,e_2,\dots)
\mid e_n\in E_n, \ t(e_n)=i(e_{n+1})\text{ for all }n\in\N\} \]
be the path space associated to $(V,E)$. 
Equipped with the relative topology 
from the product space $\prod_nE_n$, 
$X_{(V,E)}$ is compact, metrizable and zero-dimensional. 
We denote the cofinality relation on $(V,E)$ by $AF(V,E)$, 
that is, two paths are equivalent if they agree from some level on. 
We now equip $AF(V,E)$ with an AF-structure. 
Let $n\in\{0,1,2,\dots\}$. 
Then $AF_n(V,E)$ will denote 
the compact \'etale subequivalence relation of $AF(V,E)$ 
defined by the property of cofinality from level $n$ on. 
That is, if $x=(e_1,e_2,\dots,e_n,e_{n+1},\dots)$, 
$y=(f_1,f_2,\dots,f_n,f_{n+1},\dots)$ is in $X_{(V,E)}$, then 
$(x,y)\in AF_n(V,E)$ if $e_{n+1}=f_{n+1},e_{n+2}=f_{n+2},\dots$, 
and $AF_n(V,E)$ is given 
the relative topology from $X_{(V,E)}\times X_{(V,E)}$, 
thus getting a CEER structure. 
Obviously $AF_n(V,E)\subset AF_{n+1}(V,E)$, 
and we have $AF(V,E)=\bigcup_{n=0}^\infty AF_n(V,E)$. 
We give $AF(V,E)$ the inductive limit topology, i.e. 
$\displaystyle AF(V,E)=\lim_{\longrightarrow}AF_n(V,E)$. 

Let $p=(e_1,e_2,\dots,e_n)$ be a finite path 
from level 0 to some level $n$, 
and let $U(p)$ denote the \emph{cylinder set} in $X_{(V,E)}$ 
defined by 
\[ U(p)=\{x=(f_1,f_2,\dots)\in X_{(V,E)}\mid
f_1=e_1,f_2=e_2,\dots,f_n=e_n\}. \]
Then $U(p)$ is a clopen subset of $X_{(V,E)}$, 
and the collection of all cylinder sets is 
a clopen basis for $X_{(V,E)}$. 
Let $p=(e_1,e_2,\dots,e_n)$ and $q=(e'_1,e'_2,\dots,e'_n)$ be 
two finite paths from level 0 to the same level $n$, 
such that $t(e_n)=t(e'_n)$. 
Let $U(p,q)$ denote the intersection of $AF_n(V,E)$ 
with the Cartesian product $U(p)\times U(q)$. 
The collection of sets of the form $U(p,q)$ is 
a clopen basis for $AF(V,E)$. 
From this it follows immediately that 
a probability measure $\mu$ on $X_{(V,E)}$ is $AF(V,E)$-invariant 
(cf. Section 1) if and only if $\mu(U(p))=\mu(U(q))$ 
for all such cylinder sets $U(p)$ and $U(q)$. 
\bigskip

The following theorem is proved in \cite[Theorem 3.9]{GPS2}. 
However, we will explain the salient feature of the proof, 
as this will be important for arguments later in this paper. 

\begin{thm}\label{model}
Let 
$\displaystyle(R,\mathcal{T})=\lim_{\longrightarrow}(R_n,\mathcal{T}_n)$ 
be an AF-equivalence relation on $X$. 
There exists a Bratteli diagram $(V,E)$ such that 
$(R,\mathcal{T})$ is isomorphic to 
the AF-equivalence relation $AF(V,E)$ associated to $(V,E)$. 
Furthermore, $(V,E)$ is simple if and only if $(R,\mathcal{T})$ is minimal. 
\end{thm}
\begin{figure}
\begin{center}
\begin{picture}(400,170)(0,30)

\thicklines\multiput(20,50)(0,16){9}{\line(1,0){30}}\thinlines
\put(30,30){$T_1$}
\multiput(26,50)(0,16){9}{\circle*{3}}
\multiput(44,50)(0,16){9}{\circle*{3}}

\thicklines\multiput(90,50)(0,16){6}{\line(1,0){30}}\thinlines
\put(100,30){$T_2$}
\multiput(96,50)(0,16){6}{\circle*{3}}
\put(94,134){$x$}

\multiput(135,50)(8,0){3}{\circle*{2}}

\thicklines\multiput(170,50)(0,16){3}{\line(1,0){30}}\thinlines
\put(180,30){$T_m$}
\multiput(176,50)(0,16){3}{\circle*{3}}
\multiput(188,50)(0,16){3}{\circle*{3}}
\multiput(194,50)(0,16){3}{\circle*{3}}

\put(300,60){\circle*{4}}
\put(260,140){\circle*{4}}
\put(300,140){\circle*{4}}
\multiput(320,140)(10,0){3}{\circle*{2}}
\put(360,140){\circle*{4}}
\qbezier(260,132)(270,90)(292,64)
\qbezier(264,134)(290,110)(297,66)
\put(300,66){\line(0,1){68}}
\qbezier(354,134)(320,110)(304,68)
\put(306,64){\line(3,4){52}}
\qbezier(362,134)(350,92)(308,62)

\put(290,44){$w=\widetilde{T}_l$}
\put(254,146){$v_1$}
\put(256,156){$\shortparallel$}\put(254,166){$T_1$}
\put(294,146){$v_2$}
\put(296,156){$\shortparallel$}\put(294,166){$T_2$}
\put(354,146){$v_m$}
\put(356,156){$\shortparallel$}\put(354,166){$T_m$}

\end{picture}
\caption{Illustration of the maps $i,t:E\to V$}
\end{center}
\end{figure}
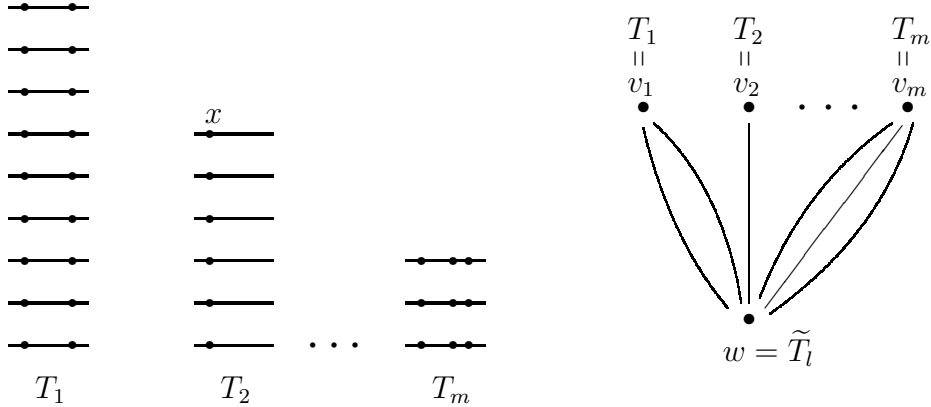%
\begin{proof}[Proof sketch]
For each $n$, choose a partition $\mathcal{P}'_n$ of $R_n$ 
such that $\bigcup_{n=0}^\infty\mathcal{P}'_n$ generates 
the topology $\mathcal{T}$ of $R$. 
(We assume that $R_0$ equals $\Delta(=\Delta_X)$, the diagonal of $X$, 
and we will freely identify $X$ with $\Delta$ whenever that is convenient.) 
Assume we have inductively obtained 
a groupoid partition $\mathcal{O}'_n$ of $R_n$ which is finer than 
both $\mathcal{P}'_n$ and $\mathcal{O}'_{n-1}\cup\{R_n\setminus R_{n-1}\}$, 
cf. Proposition \ref{groupoid}. 
Obviously $\mathcal{O}_{n+1}=\mathcal{O}'_{n+1}\cap\Delta$ is finer than 
$\mathcal{O}_n=\mathcal{O}'_n\cap\Delta$, and
$\bigcup_{i=0}^\infty\mathcal{O}_i$ generates the topology of $X$. 
Let $\mathcal{O}'_n$ be represented 
by $m$ towers $T_1,T_2,\dots,T_m$ of 
(not necessarily distinct) heights $h_1,h_2,\dots,h_m$ (cf. Figure 1). 
We let the vertex set $V_n$ at level $n$ of the Bratteli diagram $(V,E)$ 
be $V_n=\{v_1,v_2,\dots,v_m\}$, 
where $v_i$ corresponds to the tower $T_i$. 
So $[v_i]=[A]_{\mathcal{O}_n'}$ for some $A\in\mathcal{O}_n$, 
where $[\cdot]_{\mathcal{O}_n'}$ denotes the equivalence class of 
the relation $\sim_{\mathcal{O}_n'}$ on $\mathcal{O}_n'$, 
i.e. $B\sim_{\mathcal{O}_n'}C$ if there exists $V\in\mathcal{O}_n'$ 
such that $V^{-1}\cdot V=B$, $V\cdot V^{-1}=C$. 
Let $\mathcal{O}'_{n+1}$ be represented 
by $k$ towers $\widetilde{T}_1,\widetilde{T}_2,\dots,\widetilde{T}_k$ 
of heights $\tilde{h}_1,\tilde{h}_2,\dots,\tilde{h}_k$, and 
set $V_{n+1}=\{w_1,w_2,\dots,w_k\}$, 
where $w_j$ corresponds to the tower $\widetilde{T}_j$. 
Let $\mathcal{O}''_{n+1}=\{U\in\mathcal{O}'_{n+1}\mid U\subset R_n\}$. 
Clearly $\mathcal{O}''_{n+1}=\mathcal{O}'_{n+1}\cap R_n$. 
It is a simple observation that 
$\mathcal{O}''_{n+1}$ is a groupoid partition of $R_n$ 
that is finer than $\mathcal{O}'_n$, 
and that $\mathcal{O}''_{n+1}\cap\Delta=\mathcal{O}_{n+1}$. 
The set of edges $E_{n+1}$ between $V_n$ and $V_{n+1}$ 
is labelled by $\mathcal{O}_{n+1}$ modulo $\mathcal{O}''_{n+1}$, 
i.e. $E_{n+1}$ consists of 
$\sim_{\mathcal{O}''_{n+1}}$ equivalence classes 
(denoted by $[\cdot]_{\mathcal{O}''_{n+1}}$) of $\mathcal{O}_{n+1}$. 
Specifically, if $A,B\in\mathcal{O}_{n+1}$, we have 
$A\sim_{\mathcal{O}''_{n+1}}B$ 
if there exists $U\in\mathcal{O}''_{n+1}$ such that 
$U^{-1}\cdot U=A$, $U\cdot U^{-1}=B$. 
As we explained above, the vertex set $V_n$, resp. $V_{n+1}$ 
(i.e. the towers $T_1,T_2,\dots,T_m$, 
resp. $\widetilde{T}_1,\widetilde{T}_2,\dots,\widetilde{T}_k$) 
may be identified with 
the $\sim_{\mathcal{O}'_n}$ (resp. $\sim_{\mathcal{O}'_{n+1}}$) 
equivalence classes $[\cdot]_{\mathcal{O}'_n}$ 
(resp. $[\cdot]_{\mathcal{O}'_{n+1}}$) 
of $\mathcal{O}_n$ (resp. $\mathcal{O}_{n+1}$). 
Let $[A]_{\mathcal{O}''_{n+1}}\in E_{n+1}$, 
where $A\in\mathcal{O}_{n+1}$. 
Then we define $t([A]_{\mathcal{O}''_{n+1}})=[A]_{\mathcal{O}'_{n+1}}$ 
and $i([A]_{\mathcal{O}''_{n+1}})=[B]_{\mathcal{O}'_n}$, 
where $B$ is the unique element of $\mathcal{O}_n$ 
such that $A\subset B$. 
We give a more intuitive explanation of this 
by appealing to Figure 2. 
The $R_{n+1}$-equivalence class $[x]_{R_{n+1}}$ of a point $x$ in $X$ 
is marked in the tower presentation of 
the groupoid partition $\mathcal{O}'_n$ of $R_n$. 
Since $R_n\subset R_{n+1}$, 
the set $[x]_{R_{n+1}}$ breaks into 
a finite disjoint union of $R_n$-equivalence classes. 
Let $A$ be the unique clopen set in $\mathcal{O}_{n+1}$ that contains $x$, 
and let $\widetilde{T}_l$ be the unique tower 
$w\in V_{n+1}$ corresponding to $[A]_{\mathcal{O}'_{n+1}}$. 
Then $t^{-1}(w)$ consists of the edges shown in Figure 2. 
In other words, 
if $[x]_{R_{n+1}}$ contains $t_i$ distinct $R_n$-equivalence classes 
`belonging' to the tower $v_i\in V_n$, 
we connect $w$ to $v_i$ by $t_i$ edges. 
We have that $\tilde{h}_l=\sum_{i=1}^mt_ih_i$. 

In this way, we construct the Bratteli diagram $(V,E)$, and 
the map $F:X\to X_{(V,E)}$ is defined 
by $F(x)=(e_1,e_2,\dots,e_{n+1},\dots)$, 
where $e_{n+1}$ is $[A]_{\mathcal{O}''_{n+1}}$. 
It is now straightforward to show that $F$ is a homeomorphism 
and that $F\times F:R\to AF(V,E)$ is an isomorphism. 
\end{proof}

\begin{rem}\label{Ifn=0}
If $n=0$, then $\mathcal{O}'_0$ consists of one tower of height $1$ 
($\mathcal{O}'_0=\{\Delta\}$, $\mathcal{O}_0=\{X\}$). 
So $V_0$ is a one-point set, $V_0=\{v_0\}$. 
We have $V_1=\{w_1,w_2,\dots,w_k\}$, 
and $w_j$ is connected to $v_0$ by $\tilde{h}_j$ edges, 
where $\tilde{h}_j$ is the height of the tower 
corresponding to $w_j$. 
\end{rem}

\section{Transverse equivalence relations}

Let $X$ be a compact, metrizable and zero-dimensional space, 
and let $R$ and $S$ be (countable) equivalence relations on $X$. 
We let $R\vee S$ denote the (countable) equivalence relation on $X$ 
generated by $R$ and $S$. 
If $S\subset R$, i.e. $R\vee S=R$, 
we say that $S$ is a \emph{subequivalence} relation of $R$. 
If $S$ is a subequivalence relation of 
the \'etale equivalence relation $(R,\mathcal{T})$ and $S\in\mathcal{T}$ 
(i.e. $S$ is open in $R$), then 
$(S,\mathcal{T}|_S)$ (i.e. $S$ with the relative topology) is 
\'etale---a fact that is easily shown. 

For equivalence relations $R$ and $S$ on $X$, 
we define the following subset, denoted $R\times_XS$, 
of the Cartesian product $R\times S$:
\[ R\times_XS=\{((x,y),(y,z))\mid(x,y)\in R,(y,z)\in S\} \]
and we define $r,s:R\times_XS\to X$ by 
\[ r((x,y),(y,z))=x\quad\text{and}\quad s((x,y),(y,z))=z.\]
Also, we define $r\times s:R\times_XS\to X\times X$ 
by $r\times s((x,y),(y,z))=(x,z)$. 

If $R$ and $S$ have topologies, 
we give $R\times_XS$ the relative topology from $R\times S$ 
(with the product topology). 

\begin{df}[Transverse equivalence relation]\label{transverse}
Let $R$ and $S$ be \'etale equivalence relations on $X$. 
We say that $R$ and $S$ are \emph{transverse} to each other, 
denoted $R\perp S$, if the following hold: 
\begin{enumerate}
\item $R\cap S=\Delta_X(=\{(x,x)\mid x\in X\})$. 
\item There is a homeomorphism $h:R\times_XS\to S\times_XR$ 
such that $r\circ h=r$ and $s\circ h=s$. 
(Note that $h^{-1}:S\times_XR\to R\times_XS$ satisfies 
$r\circ h^{-1}=r$ and $s\circ h^{-1}=s$.) 
So, for each $(x,y)$ in $R$ and $(y,z)$ in $S$, 
there is a unique $y'$ in $X$ such that 
$(x,y')$ is in $S$, $(y',z)$ is in $R$ and 
$h((x,y),(y,z))=((x,y'),(y',z))$. 
\end{enumerate}
\end{df} 

An important class of examples is the following. 
Suppose that $R$ is an \'etale equivalence relation on $X$, 
and $\alpha:G\to\Homeo(X)$ is an action of 
the countable (or finite) group $G$ as homeomorphisms on $X$ 
such that 
\begin{enumerate}
\item $\alpha_g\times\alpha_g(R)=R$ for all $g\in G$. 
\item $\alpha_g\times\alpha_g:R\to R$ is a homeomorphism 
for all $g\in G$. 
((i) and (ii) together says that $\alpha_g$ implements an automorphism 
of $R$ for all $g\in G$.) 
\item $(x,\alpha_g(x))$ is not in $R$ for any $x$ in $X$ and 
$g\neq e$ in $G$ (in particular, the action is free). 
\end{enumerate}
Then the equivalence relation 
$R_G=\{(x,\alpha_g(x))\mid x\in X,g\in G\}$ is \'etale 
(cf. Example \ref{grp}) and transverse to $R$, 
the map $h:R\times_XR_G\to R_G\times_XR$ being defined by 
\[ h((x,y),(y,\alpha_g(y)))
=((x,\alpha_g(x)),(\alpha_g(x),\alpha_g(y))) \]
for all $(x,y)$ in $R$ and $g$ in $G$. 

\begin{lem}\label{rtimess}
Let $R$ and $S$ be transverse \'etale equivalence relations on $X$. 
The equivalence relation on $X$ generated by $R$ and $S$, 
$R\vee S$, is equal to $r\times s(R\times_XS)$, 
respectively $r\times s(S\times_XR)$. 
Furthermore, the map $r\times s:R\times_XS\to R\vee S$ 
(respectively, $r\times s:S\times_XR\to R\vee S$) is a bijection. 
(We do not need the map $h$ in Definition \ref{transverse} (ii) 
to be a homeomorphism for this proof.) 
\end{lem}
\begin{proof}
We consider the map $r\times s:R\times_XS\to R\vee S$ 
(it being obvious that the arguments we give 
apply similarly to the map $r\times s:S\times_XR\to R\vee S$, 
since by Definition \ref{transverse} (ii), 
$r\times s(R\times_XS)=r\times s(S\times_XR)$). 
Clearly $r\times s(R\times_XS)\subset R\vee S$. 
If $(x,y)\in R$, then $((x,y),(y,y))\in R\times_XS$, and so 
$(x,y)\in r\times s(R\times_XS)$. 
Hence $R\subset r\times s(R\times_XS)$. 
Likewise we show that $S\subset r\times s(R\times_XS)$. 
If $(x,z)$ equals $r\times s((x,y),(y,z))=r\times s((x,y'),(y',z))$, 
then $(x,y),(x,y')$ are in $R$, and so $(y,y')\in R$. 
Likewise, $(y,z),(y',z)$ are in $S$, and so $(y,y')\in S$. 
Hence $y=y'$, and so the map $r\times s$ is injective. 
The proof will be completed by showing that 
$K=r\times s(R\times_XS)$ is an equivalence relation on $X$. 

Clearly, $K$ is reflexive. 
To prove symmetry, assume $(x,z)\in K$. 
There exists $y\in X$ such that $(x,y)\in R$, $(y,z)\in S$, 
and so $((z,y),(y,x))\in S\times_XR$. 
Since $r\times s(S\times_XR)=r\times s(R\times_XS)$, 
we get that $(z,x)\in K$. 
Hence $K$ is symmetric. 

To prove transitivity, assume $(x,z),(z,w)\in K$. 
We must show that $(x,w)\in K$. 
There exists $y,y'\in X$ such that 
$((x,y),(y,z)),((z,y'),(y',w))\in R\times_XS$. 
This implies that $((y,z),(z,y'))\in S\times_XR$. 
Since the map $h:R\times_XS\to S\times_XR$ is a bijection, 
there exists $y''\in X$ such that $((y,y''),(y'',y'))\in R\times_XS$. 
We thus get that $(y',w),(y'',y')\in S$, 
which implies that $(y'',w)\in S$. 
Also, we have that $(x,y),(y,y'')\in R$, 
which implies that $(x,y'')\in R$. 
Hence $((x,y''),(y'',w))\in R\times_XS$, 
which implies that $(x,w)\in r\times s(R\times_XS)=K$, 
which proves transitivity. 
\end{proof}

\begin{prop}\label{transetale}
Let $(R,\mathcal{T})$ and $(S,\widetilde{\mathcal{T}})$ be 
two \'etale equivalence relations on $X$ 
which are transverse to each other. 
By the bijective map (cf. Lemma \ref{rtimess}) 
$r\times s:R\times_XS\to R\vee S$, 
which sends $((x,y),(y,z))\in R\times_XS$ to $(x,z)$, 
we transfer the topology on $R\times_XS(\subset R\times S)$ 
to $R\vee S$. 
With this topology, denoted $\mathcal{W}$, 
$R\vee S$ is an \'etale equivalence relation on $X$. 
In particular, if $R$ and $S$ are CEERs, then 
$(R\vee S,\mathcal{W})$ is a CEER. 
Furthermore, 
$\mathcal{W}$ is the unique \'etale topology on $R\vee S$, 
which extends $\mathcal{T}$ on $R$ and $\widetilde{\mathcal{T}}$ on $S$, 
i.e. the relative topologies on $R$ and $S$ are 
$\mathcal{T}$ and $\widetilde{\mathcal{T}}$, respectively. 
Also $R$ and $S$ are both open subsets of $R\vee S$. 
\end{prop}
\begin{proof}
It is easily seen that 
$R\times_XS$ is a closed subset of $R\times S$, 
and so the relative topology on $R\times_XS$, 
and consequently $\mathcal{W}$, is locally compact and metrizable. 
If $R$ and $S$ are CEERs, then clearly $R\times_XS$ is compact. 
Now $R$ (resp. $S$) is the image under $r\times s$ of 
\begin{align*}
& R\times_X\Delta=\{((x,y),(y,y))\mid(x,y)\in R\} \\
& (\text{resp. }\Delta\times_XS=\{((x,x),(x,y))\mid(x,y)\in S\}), 
\end{align*}
where $\Delta=\Delta_X$ is the diagonal of $X\times X$. 
Now $R\times_X\Delta$ (resp. $\Delta\times_XS$) is clopen in $R\times_XS$, 
since $\Delta$ is clopen in $R$ (resp. $S$), 
and so we get that $R$ (resp. $S$) is clopen in $R\vee S$. 

We now show the \'etale condition for $R\vee S$. 
Let $(x,z)\in R\vee S$, and let $y$ be the unique point in $X$ 
such that $((x,y),(y,z))\in R\times_XS$. 
A local basis at $(x,z)$ is the family of composition of graphs 
in $\mathcal{T}$ and $\widetilde{\mathcal{T}}$ (cf. Section 1) 
\[ (\widetilde{U},\tilde{\gamma},\widetilde{V})\circ(U,\gamma,V)
=(U\cap\gamma^{-1}(V\cap\widetilde{U}),\tilde{\gamma}\circ\gamma,
\tilde{\gamma}(V\cap\widetilde{U})), \]
where $(U,\gamma,V)\in\mathcal{T}$, 
$(\widetilde{U},\tilde{\gamma},\widetilde{V})\in\widetilde{\mathcal{T}}$, 
and $x\in U$, $y\in V\cap\widetilde{U}$, $z\in\widetilde{V}$, 
such that $y=\gamma(x)$, $z=\tilde{\gamma}(y)$. 
In fact, we may assume that $\widetilde{U}=V$, 
and so a local basis at $(x,z)$ is the family of graphs 
\[ \left\{(U,\tilde{\gamma}\circ\gamma,V) \ \left| \ 
\begin{matrix}
x\in U, \ z\in V, \ U\text{ and }V\text{ open in }X,\\
V=\tilde{\gamma}\circ\gamma(U), \ 
\graph(\gamma)\subset R, \ \graph(\tilde{\gamma})\subset S
\end{matrix}
\right.\right\}. \]
Clearly each $(U,\tilde{\gamma}\circ\gamma,V)$ is 
an \'etale open neighbourhood of $(x,z)$ in $\mathcal{W}$. 
We need to show that 
the product of composable pairs in $R\vee S$ is continuous, 
and also that the inverse map on $R\vee S$ is continuous. 
If we assume this has been established, 
it follows easily that 
any \'etale topology $\mathcal{E}$ on $R\vee S$ 
that extends $\mathcal{T}$ on $R$ and $\widetilde{\mathcal{T}}$ on $S$ 
has to be equal to $\mathcal{W}$. 
In fact, since $(x,y)\cdot(y,z)=(x,z)$, 
where $(x,y)\in R$, $(y,z)\in S$, 
it is easy to show that 
the graphs $(U,\tilde{\gamma}\circ\gamma,V)$ considered above is also 
a local basis for $\mathcal{E}$ at $(x,z)$. 
Thus, $\mathcal{W}$ is contained in $\mathcal{E}$. 
The other inclusion follows easily 
from the \'etaleness of $\mathcal{E}$. 

To prove that the map $(x,z)\mapsto (x,z)^{-1}=(z,x)$ is continuous 
(and hence a homeomorphism) on $R\vee S$, let $(x_n,z_n)\to(x,z)$. 
This means that there exist (unique) $y_n,y\in X$ such that 
$((x_n,y_n),(y_n,z_n))\to((x,y),(y,z))$ in $R\times_XS$. 
This implies that 
$((z_n,y_n),(y_n,x_n))\to((z,y),(y,x))$ in $S\times_XR$. 
Applying the map $h$ of Definition \ref{transverse} (ii), 
we conclude that there exist $y'_n,y'\in X$ such that 
$((z_n,y'_n),(y'_n,x_n))\to((z,y'),(y',x))$ in $R\times_XS$. 
This means that $(z_n,x_n)\to(z,x)$ in $R\vee S$, and we are done. 

To prove that the product of composable pairs 
in $(R\vee S)\times (R\vee S)$ is continuous, 
let $((x_n,y_n),(y_n,z_n))\to((x,y),(y,z))$ 
in $(R\vee S)\times (R\vee S)$. 
We want to show that $(x_n,z_n)\to(x,z)$ in $R\vee S$. 
There exist (unique) $y'_n,y''_n,y',y''\in X$ such that 
$((x_n,y'_n),(y'_n,y_n))\to((x,y'),(y',y))$ and 
$((y_n,y''_n),(y''_n,z_n))\to((y,y''),(y'',z))$ in $R\times_XS$. 
This implies that 
$((y'_n,y_n),(y_n,y''_n))\to((y',y),(y,y''))$ in $S\times_XR$. 
Using the map $h$ of Definition \ref{transverse} (ii), 
there exist $y'''_n,y'''\in X$ such that 
$((y'_n,y'''_n),(y'''_n,y''_n))\to((y',y'''),(y''',y''))$ 
in $R\times_XS$. 
So we have altogether 
$(x_n,y'_n)\to(x,y')$, $(y'_n,y'''_n)\to(y',y''')$ in $R$, 
$(y'''_n,y''_n)\to(y''',y'')$, $(y''_n,z_n)\to(y'',z)$ in $S$. 
This implies that $(x_n,y'''_n)\to(x,y''')$ in $R$, 
and $(y'''_n,z_n)\to(y''',z)$ in $S$. 
Hence $((x_n,y'''_n),(y'''_n,z_n))\to((x,y'''),(y''',z))$ 
in $R\times_XS$, and so $(x_n,z_n)\to(x,z)$ in $R\vee S$. 
\end{proof}

Henceforth, 
whenever $(R,\mathcal{T})$ and $(S,\widetilde{\mathcal{T}})$ are 
two transverse equivalence relations on $X$, 
we will give $R\vee S$ the \'etale topology $\mathcal{W}$ 
described in Proposition \ref{transetale}. 

We want to prove that if $R$ is an AF-equivalence relation on $X$, 
and $S$ is a CEER on $X$ such that $R\perp S$ 
(i.e. $R$ and $S$ are transverse), 
then $R\vee S$ is again AF. 
Furthermore, we will give an explicit description of the relation 
between the Bratteli diagram models for $R\vee S$ and $R$, respectively. 
This will be important for the proof of the absorption theorem 
in the next section. 
We shall need the following two lemmas. 

\begin{lem}\label{transCEER}
Let $(R,\mathcal{T})$ and $(S,\widetilde{\mathcal{T}})$ be 
two transverse CEERs on $X$. 
Then the following hold: 
\begin{enumerate}
\item $(x,y)\in R\Rightarrow\#([x]_S)=\#([y]_S)$. 
(In fact, $R$ does not have to be a CEER for (i) to hold.) 
\item Let $\mathcal{O}'$ be a groupoid partition of $R\vee S$ 
that is finer than the (clopen) partition 
$\{\Delta,R\setminus\Delta,S\setminus\Delta,
(R\vee S)\setminus(R\cup S)\}$. 
Then for any $U\in\mathcal{O}'$, 
there are unique elements $U_R$, $U_S$, $V_R$, $V_S$ in $\mathcal{O}'$ 
such that $U_R,V_R\subset R$, $U_S,V_S\subset S$ and 
$U=U_R\cdot U_S=V_S\cdot V_R$. 
The partitions $\mathcal{O}'|_R(=\mathcal{O}'\cap R)$ 
and $\mathcal{O}'|_S(=\mathcal{O}'\cap S)$ are groupoid partitions 
for $R$ and $S$, respectively. 
\end{enumerate}
\end{lem}
\begin{proof}
(i). 
Let $(x,y)\in R$ and let $[x]_S=\{x{=}x_1,x_2,\dots,x_n\}$, 
$[y]_S=\{y{=}y_1,y_2,\dots,y_m\}$, where $n,m\in\N$ 
(cf. Proposition \ref{ceer} (iii)). 
Let $y_i\in[y]_S$. 
Then $((x,y),(y,y_i))\in R\times_XS$. 
By Definition \ref{transverse} (ii) 
there exists a unique $x_i\in[x]_S$ such that 
$h((x,y),(y,y_i))=((x,x_i),(x_i,y_i))$. 
We will prove that 
the map $y_i\in[y]_S\mapsto x_i\in[x]_S$ is one-to-one. 
In fact, assume $y_j\in[y]_S$ and that 
$h((x,y),(y,y_j))=((x,x_i),(x_i,y_j))$. 
Since $(x_i,y_i),(x_i,y_j)\in R$, 
we get that $(y_i,y_j)\in R$. 
Likewise $(y,y_i),(y,y_j)\in S$, and so $(y_i,y_j)\in S$. 
Since $R\cap S=\Delta$, we get that $y_i=y_j$, and so $j=i$. 
We conclude that $\#([x]_S)\geq\#([y]_S)$. 
Similarly, 
by considering the map $x_j\in[x]_S\mapsto y_j\in[y]_S$ 
defined by $h((y,x),(x,x_j))=((y,y_j),(y_j,x_j))$, 
we show that $\#([y]_S)\geq\#([x]_S)$. 
Hence we get that $\#([x]_S)=\#([y]_S)$. 

(ii). 
Let $(x,z)$ be any element of $U$. 
There exists a unique $y\in X$ 
such that $(x,y)\in R$ and $(y,z)\in S$. 
Let $U_R,U_S$ be the unique elements of $\mathcal{O}'$ 
which contain $(x,y)$ and $(y,z)$, respectively. 
Since $(x,y)\in R$ and 
the partition $\mathcal{O}'$ is finer than 
the partition $\{R,(R\vee S)\setminus R\}$, 
we have $U_R\subset R$. 
Similarly, we have $(y,z)\in U_S\subset S$. 
Since $U_R\cdot U_S$ contains $(x,z)$, 
it meets $U$ and hence $U=U_R\cdot U_S$. 
The existence and uniqueness of $V_R$ and $V_S$ are shown 
in an analogous way. 

That $\mathcal{O}'|_R$ and $\mathcal{O}'|_S$ are groupoid partition 
for $R$ and $S$, respectively, is obvious. 
\end{proof}

By Lemma \ref{transCEER} 
each tower associated to a groupoid partition $\mathcal{O}'$ of $R\vee S$ 
(cf. Figure 1) is decomposed into an ``orthogonal'' array of towers 
that are associated to the groupoid partitions 
$\mathcal{O}'|_R$ and $\mathcal{O}'|_S$ of $R$ and $S$, respectively. 
In fact, 
let $T$ be one of the ``$R\vee S$-towers'' associated to $\mathcal{O}'$ 
of height $k$. 
If $x\in T$, then $[x]_{R\vee S}$ (which has cardinality $k$) is 
a disjoint union of $m$ $R$-equivalence classes, 
respectively $n$ $S$-equivalence classes, where $k=mn$. 
So $T$ is the disjoint union of $m$ $R$-towers 
associated to the groupoid partition $\mathcal{O}'|_R$ of $R$, 
and also the disjoint union of $n$ $S$-towers 
associated to the groupoid partition $\mathcal{O}'|_S$ of $S$. 
We have illustrated this in Figure 3, 
where the $R$-towers are drawn vertical and 
the $S$-towers are drawn horizontal. 
The $R\vee S$-equivalence class of $x$ is marked, and 
we indicate how a given $U\in\mathcal{O}'$ associated to the tower $T$ 
can be written as $U=U_R\cdot U_S=V_S\cdot V_R$, 
as explained in Lemma \ref{transCEER}. 

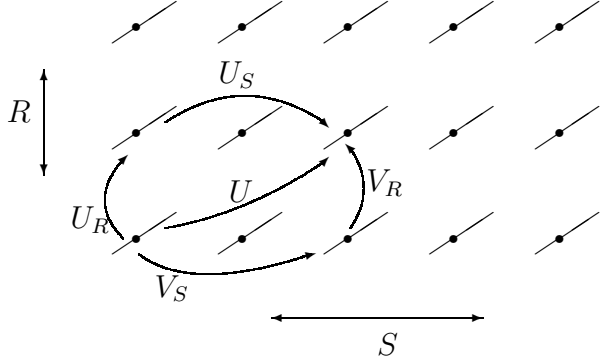
\begin{figure}
\begin{center}
\begin{picture}(400,170)(0,30)

\multiput(80,70)(40,0){5}{\line(3,2){24}}
\multiput(80,110)(40,0){5}{\line(3,2){24}}
\multiput(80,150)(40,0){5}{\line(3,2){24}}

\multiput(89,76)(40,0){5}{\circle*{3}}
\multiput(89,116)(40,0){5}{\circle*{3}}
\multiput(89,156)(40,0){5}{\circle*{3}}

\put(54,120){\vector(0,1){20}}
\put(54,120){\vector(0,-1){20}}

\put(40,120){$R$}

\put(180,46){\vector(1,0){40}}
\put(180,46){\vector(-1,0){40}}

\put(180,32){$S$}

\qbezier(100,80)(130,85)(160,105)
\put(160,105){\vector(1,1){2}}
\put(124,90){$U$}

\qbezier(100,120)(130,140)(160,120)
\put(160,120){\vector(1,-1){2}}
\put(120,135){$U_S$}

\qbezier(90,70)(110,55)(155,70)
\put(155,70){\vector(4,1){2}}
\put(96,54){$V_S$}

\qbezier(84,76)(70,90)(84,106)
\put(84,106){\vector(1,1){2}}
\put(64,80){$U_R$}

\qbezier(170,80)(180,95)(170,110)
\put(170,110){\vector(-1,1){2}}
\put(176,94){$V_R$}

\end{picture}
\caption{The decomposition of a $(R\vee S)$-tower 
into $R$-towers and $S$-towers}
\end{center}
\end{figure}

\begin{lem}\label{transCEERs}
Let 
$\displaystyle(R,\mathcal{T})=\lim_{\longrightarrow}(R_n,\mathcal{T}_n)$ 
be an AF-equivalence relation on $X$, 
where $\Delta_X{=}\Delta{=}R_0\subset R_1\subset R_2\subset\dots$ is 
an ascending sequence of CEERs on $X$. 
Let $(S,\widetilde{\mathcal{T}})$ be a CEER which is transverse to $R$, 
i.e. $R\perp S$. 
There exists an ascending sequence of CEERs $\{(R'_n,\mathcal{T}'_n)\}$, 
$\Delta{=}R'_0\subset R'_1\subset R'_2\subset\dots$, 
such that 
$\displaystyle(R,\mathcal{T})=\lim_{\longrightarrow}(R'_n,\mathcal{T}'_n)$ 
and $R'_n\perp S$ for all $n$. 
\end{lem}
\begin{proof}
Define the following subset $R'_n$ of $R_n$ by 
\[ R'_n=\left\{(x,y)\in R_n\left|
\begin{array}{l}
\forall(y,z)\in S, \\
h((x,y),(y,z))=((x,y'),(y',z))\text{ implies }
(y',z)\in R_n
\end{array}\right.\right\}, \]
where $h:R\times_XS\to S\times_XR$ is the map 
in Definition \ref{transverse}. 
By slight abuse of notation 
we may alternatively define $R'_n(\subset R_n)$ by 
\[ (x,y)\in R'_n\Longleftrightarrow
h(\{(x,y)\}\times_XS)\subset S\times_XR_n. \]
Clearly $R'_n\subset R'_{n+1}$ for all $n$. 
Also, we claim that $\bigcup_{n=0}^\infty R'_n=R$. 
In fact, let $(x,y)\in R$. 
By Proposition \ref{ceer} (iii), 
the $S$-equivalence class of $y$ is finite, 
say $[y]_S=\{y{=}y_1,y_2,y_3,\dots,y_L\}$. 
For each $1\leq l\leq L$, 
there exists a unique $x_l\in X$ 
such that $h((x,y),(y,y_l))=((x,x_l),(x_l,y_l))$, 
where $(x_l,y_l)\in R$. 
We choose $N$ sufficiently large so that 
$(x_l,y_l)\in R_N$ for all $1\leq l\leq L$. 
Since $y_1=y$, we must have $x_1=x$ 
by the properties of the map $h$. 
Hence $(x,y)\in R'_N$. 

We now prove that $R'_n$ is an equivalence relation for all $n$. 
Reflexivity is obvious from the definition of $R'_n$, 
using the fact that $h((x,x),(x,z))=((x,z),(z,z))$. 
Now let $(x,y)\in R'_n(\subset R)$. 
By Lemma \ref{transCEER} (i), we have that $\#([x]_S)=\#([y]_S)$. 
By appropriate labelling, we have for any $1\leq l\leq m$, 
\[ h((x,y),(y,y_l))=((x,x_l),(x_l,y_l)), \tag{$*$} \]
where $[x]_S=\{x{=}x_1,x_2,\dots,x_m\}$, 
$[y]_S=\{y{=}y_1,y_2,\dots,y_m\}$ and $(x_l,y_l)\in R_n$. 
Since $h((y,x),(x,x_l))=((y,y_l),(y_l,x_l))$, 
we conclude that $(y,x)\in R'_n$, and so 
$R'_n$ is symmetric. 
To prove transitivity, 
assume $(x,y),(y,z)\in R'_n$, 
and let $[x]_S=\{x{=}x_1,x_2,\dots,x_m\}$, 
$[y]_S=\{y{=}y_1,\dots,y_m\}$, $[z]_S=\{z{=}z_1,z_2,\dots,z_m\}$. 
The labelling is done as explained above (cf. ($*$)), 
so that $(x_l,y_l),(y_l,z_l)\in R_n$ for $1\leq l\leq m$. 
Since $(x_l,z_l)\in R_n$, we get that 
$h((x,z),(z,z_l))=((x,x_l),(x_l,z_l))$ for all $1\leq l\leq m$. 
This implies that $(x,z)\in R'_n$, 
thus finishing the proof that $R'_n$ is an equivalence relation. 

Note that we can use ($*$) to deduce that 
the map $h$ (or rather, its restriction), 
$h:R'_n\times _XS\to S\times_XR'_n$, is a bijection. 
In fact, let $l$ be fixed, $1\leq l\leq m$. 
We must show that $(x_l,y_l)$ in ($*$) lies in $R'_n$. 
We have that $h((y_l,x_l),(x_l,x_j))=((y_l,y_j),(y_j,x_j))$ 
for all $1\leq j\leq m$. 
Since $(y_j,x_j)\in R_n$ for $1\leq j\leq m$, 
we conclude that $(y_l,x_l)$, and hence $(x_l,y_l)$, is in $R'_n$. 

Obviously we have $R'_n\cap S=\Delta$. 
So to finish the proof, it is sufficient to show that 
$R'_n$ is a clopen subset of $R_n$. 
We claim that 
\[ R'_n\times_XS=(R_n\times_XS)\cap h^{-1}(S\times_XR_n). \tag{$**$} \]
In fact, it is clear that 
the set on the left hand side of ($**$) is contained 
in the set on the right hand side of ($**$). 
Conversely, let $((x,y),(y,z))\in (R_n\times_XS)\cap h^{-1}(S\times_XR_n)$. 
Then $(x,y)\in R_n$, and 
$h((x,y),(y,z))=((x,y'),(y',z))$ implies $(y',z)\in R_n$. 
So $(x,y)\in R'_n$, proving the other inclusion of ($**$). 

Now $R_n\times_XS$ and $S\times_XR_n$ are easily seen to be closed subsets 
of the Cartesian products $R_n\times S$ and $S\times R_n$, respectively 
(cf. Proposition \ref{ceer} (i)). 
Hence they are compact and, a fortiori, closed subsets 
of $R\times_XS$ and $S\times_XR$, respectively. 
Since $R_n$ is open in $R$, 
it follows easily that $R_n\times_XS$ and $S\times_XR_n$ are open subsets 
of $R\times_XS$ and $S\times_XR$, respectively. 
From ($**$) we conclude that 
$R'_n\times_XS$ is a compact and open subset of $R\times_XS$. 
Now $R'_n=\pi_1(R'_n\times_XS)$, 
where $\pi_1:R\times S\to R$ is the projection map. 
Since $\pi_1$ is a continuous and open map, 
we conclude that $R'_n$ is compact and open in $R$, 
and hence in $R_n$. 
This completes the proof. 
\end{proof}

\begin{prop}\label{transBratteli}
Let 
$\displaystyle(R,\mathcal{T})=\lim_{\longrightarrow}(R_n,\mathcal{T}_n)$ 
be an AF-equivalence relation on $X$, 
and let $S$ be a CEER on $X$ such that $R\perp S$, 
i.e. $R$ and $S$ are transverse to each other. 
Then $R\vee S$ is an AF-equivalence relation. 

Furthermore, there exist Bratteli diagrams $(V,E)$ and $(V',E')$ such that 
\[ R\cong AF(V,E),\quad R'=R\vee S\cong AF(V',E') \]
and so that $t:E_1\to V_1$ is injective and $S\cong AF_1(V',E')$. 
Moreover, there are surjective maps (respecting gradings) 
$q_V:V\to V'$ and $q_E:E\to E'$ 
such that 
\begin{enumerate}
\item $i(q_E(e))=q_V(i(e))$, $t(q_E(e))=q_V(t(e))$ for $e$ in $E$. 
\item For each $v$ in $V$, 
$q_E:i^{-1}(\{v\})\to i^{-1}(\{q_V(v)\})$ is a bijection. 
\item For each $v\in V_n$ and $n\geq2$, 
$q_E:t^{-1}(\{v\})\to t^{-1}(\{q_V(v)\})$ is a bijection. 
\end{enumerate}
The map $H:X_{(V,E)}\to X_{(V',E')}$ defined by 
\[ x=(e_1,e_2,e_3,\dots)\mapsto H(x)=(q_E(e_1),q_E(e_2),q_E(e_3),\dots) \]
is a homeomorphism, 
and $H$ implements an embedding of $AF(V,E)$ into $AF(V',E')$ 
whose image is transverse to $AF_1(V',E')\cong S$. 
Moreover, we have 
\[ AF_1(V',E')\vee(H\times H)(AF(V,E))
\cong AF(V',E'). \]
(Recall the notation and terminology we introduced in Section 2. 
In particular, $V=V_0\cup V_1\cup V_2\cup\dots$, $E=E_1\cup E_2\cup\dots$, 
$V'=V'_0\cup V'_1\cup V'_2\cup\dots$, $E'=E'_1\cup E'_2\cup\dots$.)
\end{prop}
\begin{proof}
By Lemma \ref{transCEERs}, we may assume that $R_n\perp S$ for all $n$. 
Now let $(\Delta_X=)\Delta=R_0=R_1\subset R_2\subset R_3\subset\dots$ and 
\[ \Delta{=}R'_0\subset R'_1{=}R_1\vee S\subset 
R'_2{=}R_2\vee S\subset R'_3{=}R_3\vee S\subset\dots\subset R'=R\vee S
=\bigcup_{n=0}^{\infty}R'_n. \]
Applying Proposition \ref{transetale} we can conclude that 
$R\vee S$ is an AF-equivalence relation. 

Now let $\mathcal{P}'_n$ be a groupoid partition of $R'_n=R_n\vee S$ 
which is finer than 
both $\mathcal{P}'_{n-1}\cup\{R'_n\setminus R'_{n-1}\}$ and 
$\{\Delta,R_n\setminus\Delta,S\setminus\Delta,R'_n\setminus(R_n\cup S)\}$, 
where we set $\mathcal{P}'_0=\{\Delta\}$. 
We require that 
$\bigcup_{n=0}^\infty\mathcal{P}'_n$ generates the topology of $R\vee S$. 
(This can be achieved 
by successive applications of Proposition \ref{groupoid}.) 
By Lemma \ref{transCEER} we get that 
$\mathcal{O}'_n=\mathcal{P}'_n|_{R_n}(=\mathcal{P}'_n\cap R_n)$ is 
a groupoid partition of $R_n$ for $n\geq0$, 
such that $\mathcal{O}'_n$ is finer than 
both $\mathcal{O}'_{n-1}\cup\{R_n\setminus R_{n-1}\}$ and 
$\{\Delta,R_n\setminus\Delta\}$. 
Furthermore, 
$\bigcup_{n=0}^\infty\mathcal{O}'_n$ generates the topology of $R$. 
Similarly, $\mathcal{Q}'_n=\mathcal{P}'_n|_S(=\mathcal{P}'_n\cap S)$ is 
a groupoid partition of $S$ for $n\geq1$, 
such that $\mathcal{Q}'_n$ is finer than 
both $\mathcal{Q}'_{n-1}$ and $\{\Delta,S\setminus\Delta\}$, and 
$\bigcup_{n=0}^\infty\mathcal{Q}'_n$ generates the topology of $S$. 
(We set $\mathcal{Q}'_0=\{\Delta\}$.) 
By Lemma \ref{transCEER}, 
every $U\in\mathcal{P}'_n$ can be uniquely written as $U=U_R\cdot U_S$, 
where $U_R\in\mathcal{O}'_n$ and $U_S\in\mathcal{Q}'_n$, 
and we may suggestively write 
$\mathcal{P}_n'=\mathcal{O}_n'\cdot\mathcal{Q}_n'$. 
For each $n\geq0$ we have that 
\[ (\mathcal{P}'_n\cap\Delta=)\mathcal{P}'_n|_\Delta
=\mathcal{O}'_n|_\Delta=\mathcal{Q}'_n|_\Delta=\mathcal{P}_n \] 
is a clopen partition of $X$, with 
\[ \Delta=\mathcal{P}_0\prec\mathcal{P}_1\prec\mathcal{P}_2\prec\dots \]
and $\bigcup_{n=0}^\infty\mathcal{P}_n$ being a basis for $X$. 

Combining all this---following the description 
given in the proof sketch of Theorem \ref{model}---
we construct the Bratteli diagrams $(V',E')$ and $(V,E)$, 
so that $R'=R\vee S\cong AF(V',E')$ and $R\cong AF(V,E)$, respectively, 
and such that the conditions stated in the proposition are satisfied. 
For brevity we will omit some of the details, 
which are routine verifications, 
and focus on the main ingredients of the proof. 
(We will use the same notation 
that we used in the proof sketch of Theorem \ref{model}.)

First we observe that $S\cong AF_1(V',E')$. 
In fact, since $\mathcal{P}_0'=\mathcal{Q}_0'=\{\Delta\}$, 
and $\mathcal{P}_1'=\mathcal{Q}_1'$ is a groupoid partition of 
$\mathcal{R}_1'=\Delta\vee S=S$, 
the edge set $E_1'$ is related in an obvious way to the towers 
$\{T_1,T_2,\dots,T_m\}$ associated to $\mathcal{Q}_1'$ 
(cf. Figures 1 \& 2 and Remark \ref{Ifn=0}). 
The towers associated to $\mathcal{Q}_n'$, $n>1$, are obtained 
by subdividing (vertically) the towers $\{T_1,T_2,\dots,T_m\}$, 
and from this it is easily seen that 
$S$ is isomorphic to $AF_1(V',E')$. 
Since $R_0=R_1=\{\Delta\}$, and $\mathcal{O}_0'=\{\Delta\}$, 
$\mathcal{O}_1'=\mathcal{O}_1'|_\Delta=\mathcal{P}_1$, 
we deduce that $t:E_1\to V_1$ is injective and that 
the obviously defined maps, 
$q_V:V_0\to V_0'$, $q_V:V_1\to V_1'$, $q_E:E_1\to E_1'$, satisfy 
condition (i) for $e\in E_1$, and (ii) for $v\in V_0$. 
(For instance, if $e\in E_1$ corresponds to $A\in\mathcal{P}_1$, 
then it is mapped to $e'\in E_1'$, 
which corresponds to the ``floor'' $A$, 
that lies in one of the towers $\{T_1,T_2,\dots,T_m\}$.) 

Let $n\geq1$. 
Assume that we have defined $q_V:V_i\to V_i'$, $i=0,1,\dots,n$, 
and $q_E:E_j\to E_j'$, $j=1,2,\dots,n$, such that 
(i) is true for $e\in E_j$, $j=1,2,\dots,n$, 
and (ii) is true for $v\in V_i$, $i=0,1,\dots,n-1$. 
Assume also that 
$H((e_1,e_2,\dots,e_n))=(q_E(e_1),q_E(e_2),\dots,q_E(e_n))$ is 
a bijection between finite paths of length $n$ 
from the top vertices of $(V,E)$ and $(V',E')$, respectively. 

We now define $q_E:E_{n+1}\to E_{n+1}'$. 
Let $e=[A]_{\mathcal{O}_{n+1}''}\in E_{n+1}$ 
for some $A\in\mathcal{P}_{n+1}
(=\mathcal{P}_{n+1}'|_\Delta=\mathcal{O}_{n+1}'|_\Delta)$, 
where $\mathcal{O}_{n+1}''=\{U\in\mathcal{O}_{n+1}'\mid U\subset R_n\}$. 
We define $q_E(e)=[A]_{\mathcal{P}_{n+1}''}\in E_{n+1}'$, 
where $\mathcal{P}_{n+1}''=\{W\in\mathcal{P}_{n+1}'\mid W\subset R_n'\}$. 
Since $\mathcal{O}_{n+1}''\subset\mathcal{P}_{n+1}''$, 
we get that $q_E:E_{n+1}\to E_{n+1}'$ is well-defined and surjective. 
We define the map $q_V:V_{n+1}\to V_{n+1}'$ by 
$q_V(v)=[B]_{\mathcal{P}_{n+1}'}\in V_{n+1}'$, 
where $v=[B]_{\mathcal{O}_{n+1}'}\in V_{n+1}$ 
for some $B\in\mathcal{P}_{n+1}$. 
Since $\mathcal{O}_{n+1}'\subset\mathcal{P}_{n+1}'$, 
we get that the map $q_V$ is well-defined and surjective. 
It is easy to see that (i) is satisfied for all $e\in E_{n+1}$. 

We show that the maps defined in (ii) and (iii) are surjective. 
Let $v=[B]_{\mathcal{O}_n'}\in V_n$, 
where $B\in\mathcal{P}_n$ (resp. $w=[\widetilde{B}]_{\mathcal{O}_{n+1}'}
\in V_{n+1}$, where $\widetilde{B}\in\mathcal{P}_{n+1}$). 
So $q_V(v)=[B]_{\mathcal{P}_n'}\in V_n'$ 
(resp. $q_V(w)=[\widetilde{B}]_{\mathcal{P}_{n+1}'}\in V_{n+1}'$). 
Let $e'=[A]_{\mathcal{P}_{n+1}''}\in i^{-1}(\{q_V(v)\})\subset E_{n+1}'$, 
where $A\in\mathcal{P}_{n+1}$ 
(resp. $f'=[\widetilde{A}]_{\mathcal{P}_{n+1}''}
\in t^{-1}(\{q_V(w)\})\subset E_{n+1}'$, 
where $\widetilde{A}\in\mathcal{P}_{n+1}$). 
This means that $A\sim_{\mathcal{P}_n''}A'\subset B$, 
for some $A'\in\mathcal{P}_{n+1}$ 
(resp. $\widetilde{A}\sim_{\mathcal{P}_{n+1}'}\widetilde{B}$, 
or, equivalently, $[\widetilde{A}]_{\mathcal{P}_{n+1}'}
=[\widetilde{B}]_{\mathcal{P}_{n+1}'}$).  
Let $e=[A']_{\mathcal{O}_{n+1}''}\in E_{n+1}$ 
(resp. $f=[\widetilde{A}]_{\mathcal{O}_{n+1}''}\in E_{n+1}$). 
Then clearly $i(e)=v$ and $q_E(e)=e'$, 
proving that $q_E:i^{-1}(\{v\})\to i^{-1}(\{q_V(v)\})$ is surjective. 
Clearly $q_E(f)=f'$. 
Also, since $n+1\geq2$, we have that 
$\mathcal{P}_{n+1}''|_{R_1'}\succ\mathcal{Q}_1'
(=\mathcal{P}_1'=\mathcal{P}_1'')$, and so 
we may choose $\widetilde{A}$ such that 
$\widetilde{A}$ and $\widetilde{B}$ are contained 
in the same set $C\in\mathcal{P}_1$. 
This fact, 
together with $\widetilde{A}\sim_{\mathcal{P}_{n+1}'}\widetilde{B}$, 
implies that $\widetilde{A}\sim_{\mathcal{O}_{n+1}'}\widetilde{B}$. 
In fact, if $U=U_R\cdot U_S\in\mathcal{P}_{n+1}'$ 
such that $U^{-1}\cdot U=\widetilde{A}$, $U\cdot U^{-1}=\widetilde{B}$, 
with $U_R\in\mathcal{O}_{n+1}'$, $U_S\in\mathcal{Q}_{n+1}'$, 
then $U=U_R$, since $U_S$ must be the identity map on $\widetilde{B}$. 
Hence we get that 
$[\widetilde{A}]_{\mathcal{O}_{n+1}'}
=[\widetilde{B}]_{\mathcal{O}_{n+1}'}=w$, and so we have proved that 
the map $q_E:t^{-1}(\{w\})\to t^{-1}(\{q_V(w)\})$ is surjective. 

We prove that (ii) holds for $v\in V_n$. 
Since $q_E:i^{-1}(\{v\})\to i^{-1}(\{q_V(v)\})$ is surjective, 
we need to prove injectivity. 
So let $e_1,e_2\in i^{-1}(\{v\})$, and assume $q_E(e_1)=q_E(e_2)$. 
We must show that $e_1=e_2$. 
Now $e_1=[A]_{\mathcal{O}_{n+1}''}$, $e_2=[B]_{\mathcal{O}_{n+1}''}$ 
for some $A,B\in\mathcal{P}_{n+1}$. 
Since $q_E(e_1)=q_E(e_2)$, 
we have that $[A]_{\mathcal{P}_{n+1}''}=[B]_{\mathcal{P}_{n+1}''}$. 
Hence there exists $U\in\mathcal{P}_{n+1}'$ such that 
$U^{-1}\cdot U=A$, $U\cdot U^{-1}=B$, and $U\subset R_n\vee S=R_n'$. 
Since $e_1,e_2\in i^{-1}(\{v\})$, 
we have that $A\subset A_1\in\mathcal{P}_n$, $B\subset B_1\in\mathcal{P}_n$, 
such that $[A_1]_{\mathcal{O}_n'}=[B_1]_{\mathcal{O}_n'}$; 
that is, there exists $U_1\in\mathcal{O}_n'\subset\mathcal{P}_n'$, 
such that $U_1^{-1}\cdot U_1=A_1$, $U_1\cdot U_1^{-1}=B_1$. 
Since $\mathcal{P}_{n+1}'|_{R_n'}$ is finer than $\mathcal{P}_n'$, 
we must have $U\subset U_1(\subset R_n)$. 
This means that $U\in\mathcal{O}_{n+1}''$, 
and so $e_1=[A]_{\mathcal{O}_{n+1}''}=[B]_{\mathcal{O}_{n+1}''}=e_2$. 

In a similar way we prove that (iii) holds. 
In fact, let $e_1,e_2\in t^{-1}(\{v\})$, $v\in V_{n+1}$, such that 
$q_E(e_1)=q_E(e_2)$. 
We must show that $e_1=e_2$. 
Again we write 
$e_1=[A]_{\mathcal{O}_{n+1}''}$, $e_2=[B]_{\mathcal{O}_{n+1}''}$ 
for some $A,B\in\mathcal{P}_{n+1}$. 
Since $e_1,e_2\in t^{-1}(\{v\})$, 
we have $[A]_{\mathcal{O}_{n+1}'}=[B]_{\mathcal{O}_{n+1}'}$, that is, 
there exists $U_1\in\mathcal{O}_{n+1}'\subset\mathcal{P}_{n+1}'$ 
such that 
$U_1^{-1}\cdot U_1=A$, $U_1\cdot U_1^{-1}=B$. 
Since $q_E(e_1)=q_E(e_2)$, 
we have $[A]_{\mathcal{P}_{n+1}''}=[B]_{\mathcal{P}_{n+1}''}$, that is, 
there exists $U\in\mathcal{P}_{n+1}'$ such that 
$U^{-1}\cdot U=A$, $U\cdot U^{-1}=B$, and $U\subset R_n\vee S=R_n'$. 
This implies that $U=U_1$. 
Since $U_1\subset R_{n+1}$, 
$U\subset R_n\vee S$ and $S\cap R_{n+1}=\Delta$, 
we must have $U_1\subset R_n$. 
Hence $U_1\in\mathcal{O}_{n+1}''$, and so $e_1=e_2$. 

It is now straightforward to verify that 
$H:X_{(V,E)}\to X_{(V',E')}$ is a homeomorphism 
that implements an isomorphism 
between $AF(V,E)$ and its image $H\times H(AF(V,E))$ in $AF(V',E')$. 
The last assertion of the proposition is now routinely verified. 
\end{proof}

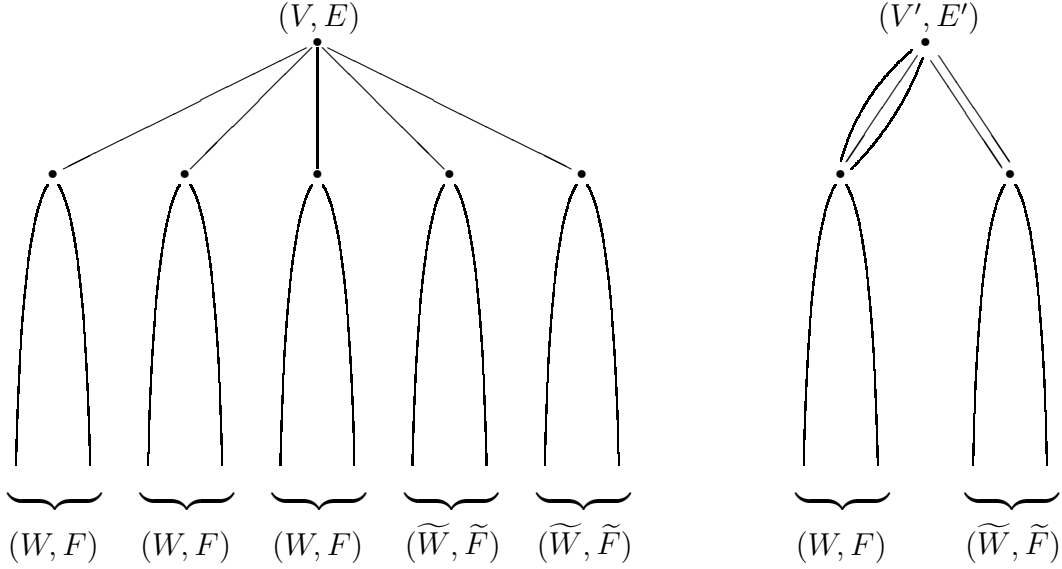
\begin{figure}
\begin{center}
\begin{picture}(400,220)(0,25)

\put(105,226){$(V,E)$}
\put(120,220){\circle*{3}}
\multiput(20,170)(50,0){5}{\circle*{3}}
\put(116,218){\line(-2,-1){92}}
\put(118,218){\line(-1,-1){46}}
\put(120,218){\line(0,-1){46}}
\put(122,218){\line(1,-1){46}}
\put(124,218){\line(2,-1){92}}
\qbezier(18,166)(8,150)(6,60)
\qbezier(68,166)(58,150)(56,60)
\qbezier(118,166)(108,150)(106,60)
\qbezier(168,166)(158,150)(156,60)
\qbezier(218,166)(208,150)(206,60)
\qbezier(22,166)(32,150)(34,60)
\qbezier(72,166)(82,150)(84,60)
\qbezier(122,166)(132,150)(134,60)
\qbezier(172,166)(182,150)(184,60)
\qbezier(222,166)(232,150)(234,60)
\put(3,50){$\underbrace{\qquad\quad}$}
\put(3,26){$(W,F)$}
\put(53,50){$\underbrace{\qquad\quad}$}
\put(53,26){$(W,F)$}
\put(103,50){$\underbrace{\qquad\quad}$}
\put(103,26){$(W,F)$}
\put(153,50){$\underbrace{\qquad\quad}$}
\put(153,26){$(\widetilde{W},\widetilde{F})$}
\put(203,50){$\underbrace{\qquad\quad}$}
\put(203,26){$(\widetilde{W},\widetilde{F})$}

\put(332,226){$(V',E')$}
\put(350,220){\circle*{3}}
\put(318,170){\circle*{3}}
\put(382,170){\circle*{3}}
\qbezier(345,217)(325,200)(318,175)
\put(347,215){\line(-2,-3){27}}
\qbezier(349,214)(340,190)(322,172)
\put(352,213){\line(2,-3){27}}
\put(355,215){\line(2,-3){27}}
\qbezier(316,166)(306,150)(304,60)
\qbezier(320,166)(330,150)(332,60)
\qbezier(380,166)(370,150)(368,60)
\qbezier(384,166)(394,150)(396,60)
\put(301,50){$\underbrace{\qquad\quad}$}
\put(301,26){$(W,F)$}
\put(365,50){$\underbrace{\qquad\quad}$}
\put(365,26){$(\widetilde{W},\widetilde{F})$}

\end{picture}
\caption{Illustrating the content of Proposition \ref{transBratteli}}
\end{center}
\end{figure}

\begin{rem}
It is helpful to illustrate 
by a figure what Proposition \ref{transBratteli} says, 
and which at the same time gives the heuristics of the proof. 
Furthermore, the illustration will be useful for 
easier comprehending the proof of the main theorem in the next section. 
In Figure 4 we have drawn the diagrams of $(V,E)$ and $(V',E')$. 
The replicate diagrams $(W,F)$ 
(resp. $(\widetilde{W},\widetilde{F})$) are also drawn, 
and $S$, the CEER transverse to $R=AF(V,E)$, 
as well as the maps $q_V$ and $q_E$, should be obvious from the figure. 
One sees that what is going on is a ``glueing'' process 
in the sense that distinct $R$-equivalence classes are ``glued'' 
together by $S$ to form $R'$-equivalence classes. 
\end{rem}

\section{The absorption theorem }

\begin{df}[$R$-\'etale and $R$-thin sets]
Let $(R,\mathcal{T})$ be an \'etale equivalence relation 
on the compact, metrizable and zero-dimensional space $X$. 
Let $Y$ be a closed subset of $X$. 
We say that 
$Y$ is \emph{$R$-\'etale} 
if the \emph{restriction} $R\cap(Y\times Y)$ of $R$ to $Y$, 
denoted by $R|_Y$, is an \'etale equivalence relation 
in the relative topology. 

We say that $Y$ is $R$-thin 
if $\mu(Y)=0$ for all $R$-invariant probability measures $\mu$ 
(cf. Section 1). 
\end{df}

\begin{rem}
It can be proved (cf. Theorem 3.11 of \cite{GPS2}) that 
if $(R,\mathcal{T})$ is an AF-equivalence relation and 
$Y$ is $R$-\'etale, then $R|_Y$ is an AF-equivalence relation on $Y$. 
Furthermore, there exist a Bratteli diagram $(V,E)$ and 
a Bratteli subdiagram $(W,F)$ 
such that $R\cong AF(V,E)$, $R|_Y\cong AF(W,F)$. 
(By a Bratteli subdiagram $(W,F)$ of $(V,E)$ 
we mean that $F$ is a subset of $E$ such that 
$i(F)=\{v_0\}\cup t(F)$, that is, 
$(W,F)$ is a Bratteli diagram, where $W=i(F)$ and $W_0=V_0$. 
We say that $F$ induces an (edge) subdiagram of $(V,E)$.
Observe that, in general, if $(W,F)$ is a subdiagram of $(V,E)$, 
then $R|_Y\cong AF(W,F)$, where $R=AF(V,E)$ and $Y=X_{(W,F)}$.) 
\end{rem}

We state a result from \cite{GPS2} 
that will be crucial in proving the absorption theorem. 

\begin{thm}[Lemma 4.15 of \cite{GPS2}]\label{glue}
Let $(R_1,\mathcal{T}_1)$ and $(R_2,\mathcal{T}_2)$ be 
two minimal AF-equivalence relations 
on the Cantor sets $X_1$ and $X_2$, respectively. 
Let $Y_i$ be a closed $R_i$-\'etale and $R_i$-thin subset of $X_i$, 
$i=1,2$. 
Assume 
\begin{enumerate}
\item $R_1\cong R_2$. 
\item There exists a homeomorphism $\alpha:Y_1\to Y_2$ 
which implements an isomorphism between $R_1|_{Y_1}$ and $R_2|_{Y_2}$. 
\end{enumerate}
Then there exists a homeomorphism $\tilde{\alpha}:X_1\to X_2$ 
which implements an isomorphism between $R_1$ and $R_2$, 
such that $\tilde{\alpha}|_{Y_1}=\alpha$, 
i.e. $\tilde{\alpha}$ is an extension of $\alpha$. 
\end{thm}

The following lemma is a technical result---easily proved---that 
we shall need for the proof of the absorption theorem. 
In the sequel 
we will use the term ``microscoping'' (of a Bratteli diagram) 
in the restricted sense called ``symbol splitting'' 
in \cite[Section 3]{GPS1}. 
(Microscoping is a converse operation to that of ``telescoping''.) 

\begin{lem}\label{enoughrooms}
Let $(R,\mathcal{T})$ be a minimal AF-equivalence relation 
on the Cantor set $X$. 
Let $\{a_n\}_{n=1}^\infty$ and $\{b_n\}_{n=1}^\infty$ be 
sequences of natural numbers. 
Then there exists a (simple) Bratteli diagram $(V,E)$ such that 
$R\cong AF(V,E)$, and for each $n\geq1$ 
\begin{enumerate}
\item $\#(V_n)\geq a_n$. 
\item For all $v\in V_{n-1}$ and all $w\in V_n$, 
$\#(\{e\in E_n|i(e){=}v,t(e){=}w\})\geq b_n$. 
\end{enumerate}
\end{lem}
\begin{proof}
Let $R\cong AF(W,F)$ for some (simple) Bratteli diagram $(W,F)$. 
By a finite number of telescopings and microscopings of 
the diagram $(W,F)$, we get $(V,E)$ with the desired properties 
(cf. \cite[Section 3]{GPS1}). 
\end{proof}

\begin{rem}\label{tele&micro}
We make the general remark that 
telescoping or microscoping a Bratteli diagram do not alter 
any of the essential properties attached to the diagram, 
like the associated path space and the AF-equivalence relation. 
In fact, there is a natural map 
between the path spaces of the original and the new Bratteli diagrams 
which implements an isomorphism 
between the AF-equivalence relations associated to the two diagrams. 
However, the versatility 
that these operations (i.e. telescopings and microscopings) give us 
in changing a given Bratteli diagram into one which is more suitable 
for our purpose---like having enough ``room'' to admit 
appropriate subdiagrams---is very helpful. 
This will be utilized extensively in the proof of the absorption theorem. 
Note that 
a subdiagram $(W,F)$ of a Bratteli diagram $(V,E)$ is being 
telescoped or microscoped (in an obvious way) simultaneously 
as these operations are applied to $(V,E)$. 
For this reason we will sometimes, when it is convenient, 
retain the old notation for the new diagrams, 
and this should not cause any confusion. 
\end{rem}

We can now state and prove the main result of this paper. 

\begin{thm}[The absorption theorem]\label{main}
Let $R=(R,\mathcal{T})$ be a minimal AF-equivalence relation 
on the Cantor set $X$, 
and let $Y$ be a closed $R$-\'etale and $R$-thin subset of $X$. 
Let $K=(K,\mathcal{S})$ be 
a compact \'etale equivalence relation on $Y$. 
Assume $K\perp R|_Y$, i.e. $K$ is transverse to $R|_Y$. 

Then there is a homeomorphism $h:X\to X$ such that 
\begin{enumerate}
\item $h\times h(R\vee K)=R$, 
where $R\vee K$ is the equivalence relation on $X$ 
generated by $R$ and $K$. 
In other words, $R\vee K$ is orbit equivalent to $R$, 
and, in particular, $R\vee K$ is affable. 
\item $h(Y)$ is $R$-\'etale and $R$-thin. 
\item $h|_Y\times h|_Y:(R|_Y)\vee K\to R|_{h(Y)}$ is a homeomorphism. 
\end{enumerate}
\end{thm}
\begin{proof}
Roughly speaking, 
the idea of the proof is 
to define an (open) AF-subequivalence relation $\overline{R}$ of $R$, 
thereby setting the stage for applying Theorem \ref{glue} 
(with $R_1=R_2=\overline{R}$), 
and in the process ``absorbing'' $Y$ (and thereby $K$) 
so that $R\vee K$ becomes $R$. 
To define $\overline{R}$ we will manipulate Bratteli diagrams, 
applying Lemma \ref{enoughrooms} together with 
Proposition \ref{transBratteli}. 
The proof is rather technical, 
and to facilitate the understanding and get the main idea of the proof 
it will be helpful to have a very special, but telling, example in mind. 
We refer to Remark \ref{simplecase} for details on this. 

Let $(W,F)$ and $(W',F')$ be two (fixed) Bratteli diagrams such that 
\[ R|_Y\cong AF(W,F), \quad 
R|_Y\vee K\cong AF(W',F') \]
and let 
\[ q_W:W\to W', \quad q_F:F\to F', \quad H:X_{(W,F)}\to X_{(W',F')} \]
be maps satisfying the conditions of 
Proposition \ref{transBratteli}. 
By Theorem 3.11 of \cite{GPS2} 
there exist a (simple) Bratteli diagram $(V,E)$ and 
a subdiagram $(\widetilde{W},\widetilde{F})$ such that 
we may assume from the start that 
$X=X_{(V,E)}$, $Y=X_{(\widetilde{W},\widetilde{F})}$, 
$R=AF(V,E)$, $R|_Y=AF(\widetilde{W},\widetilde{F})$. 
By a finite number of telescopings and microscopings applied to $(V,E)$, 
and hence to $(\widetilde{W},\widetilde{F})$, we may assume that 
for all $n\geq1$, $\#(\widetilde{W}_n)\leq\frac{1}{2}\#(V_n)$, 
and for $v_0\in V_0=\widetilde{W}_0$, $w\in\widetilde{W}_n$, 
\[ \#(\{\text{paths in }X_{(\widetilde{W},\widetilde{F})}
\text{ from }v_0\text{ to }w\})
\leq\frac{1}{2}\#(\{\text{paths in }X_{(V,E)}
\text{ from }v_0\text{ to }w\}), \]
(cf. Lemma 4.12 of \cite{GPS2}). 
Furthermore, by applying Lemma \ref{enoughrooms} 
we may assume the following holds for all $n\geq1$: 
\begin{enumerate}
\item[(1)] $\displaystyle \#(V_n)\geq\#(\widetilde{W}_n)+1
+\sum_{k=1}^{n-1}\#(W_k)$ 
\item[(2)] For $v\in V_{n-1}$, $w\in V_n$, we have 
\[ \#(\{e\in E_n\mid i(e)=v, \ t(e)=w\}) \\
\geq2\sum_{k=1}^{n-1}\#(F_k) \ 
\left(\geq2\sum_{k=1}^{n-1}\#(W_k)\right). \]
\end{enumerate}
(Note that the inequalities in (1) and (2) hold, 
if we substitute $W'_k$ for $W_k$ and $F'_k$ for $F_k$, 
cf. Proposition \ref{transBratteli}.) 
Using this we are going to construct a subdiagram of $(V,E)$ 
which will consist of countable replicas of $(W',F')$. 
This subdiagram will be instrumental 
in defining the AF-subequivalence relation $\overline{R}$ of $R$ 
that we alluded to above. 
We first choose $x_\infty=(e_1,e_2,\dots)\in X_{(V,E)}$ such that 
$t(e_n)\notin \widetilde{W}_n$ for all $n$. 
At level $n$ we choose a replica of $(W',F')$ 
``emanating'' from the vertex $t(e_n)$. 
More precisely, we let $(W',F')_n$ denote the subdiagram of $(V,E)$ 
consisting of the edges $e_1,e_2,\dots,e_n$ and then 
(a replica of) $(W',F')$ starting at the vertex $t(e_n)$. 
This can be done by (1) and (2). 
Also, by (1) we have enough ``room'' so that 
we may choose the various $(W',F')_n$'s such that 
at each level $n\geq1$, the vertex sets belonging to 
$(\widetilde{W},\widetilde{F}),(W',F')_1,(W',F')_2,\dots,(W',F')_n$ 
are pairwise disjoint. 
In Figure 5, we have illustrated this. 

\begin{figure}
\begin{center}
\begin{picture}(340,240)

\put(160,210){\circle*{5}}
\qbezier(60,30)(100,170)(160,210)
\qbezier(260,30)(220,170)(160,210)
\put(145,220){$(V,E)$}

\put(121,170){\line(1,0){78}}
\put(97,130){\line(1,0){126}}
\put(80,90){\line(1,0){160}}
\put(66,50){\line(1,0){188}}

\put(199,170){\circle*{4}}
\put(223,130){\circle*{4}}
\put(240,90){\circle*{4}}
\put(254,50){\circle*{4}}

\put(186,190){\small{$e_1$}}
\put(216,150){\small{$e_2$}}
\put(235,110){\small{$e_3$}}
\put(252,70){\small{$e_4$}}

\qbezier(160,210)(100,140)(80,30)
\qbezier(160,210)(120,160)(95,30)
\put(75,26){$\underbrace{\qquad}$}
\put(66,8){\footnotesize{$(\widetilde{W},\widetilde{F})$}}

\qbezier(199,170)(160,140)(130,30)
\qbezier(199,170)(170,140)(140,30)
\put(125,26){$\underbrace{\quad}$}
\put(110,8){\footnotesize{$(W',F')_1$}}

\qbezier(223,130)(196,110)(180,30)
\qbezier(223,130)(206,110)(190,30)
\put(175,26){$\underbrace{\quad}$}
\put(162,8){\footnotesize{$(W',F')_2$}}

\qbezier(240,90)(232,70)(224,30)
\qbezier(240,90)(236,60)(234,30)
\put(216,26){$\underbrace{\quad}$}
\put(214,8){\footnotesize{$(W',F')_3$}}

\qbezier(254,50)(250,40)(246,30)
\qbezier(254,50)(252,40)(250,30)

\end{picture}
\caption{Constructing subdiagrams of $(V,E)$}
\end{center}
\end{figure}
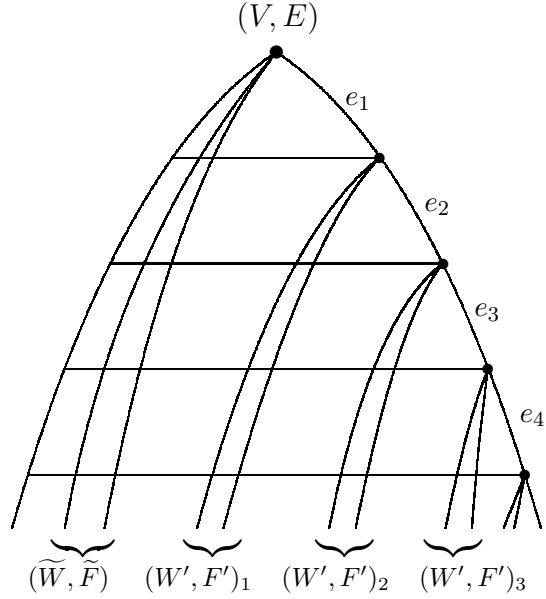

By (2) it is easily seen that 
the subdiagram $(L,G)$ of $(V,E)$ 
whose edge set consists of the edge set $\widetilde{F}$ and 
the union of the edge sets belonging to $(W',F')_n$, $n\geq1$, 
is a \emph{thin subdiagram} of $(V,E)$, 
i.e. $\mu(X_{(L,G)})=0$ for all $R$-invariant probability measures $\mu$ 
(cf. Remark \ref{tele&micro}). 
Likewise, the subdiagram $(L',G')$ 
whose edge set consists of the union of the edge sets 
belonging to $(W',F')_n$, $n\geq1$, is a thin subdiagram. 

We now construct a new Bratteli diagram $(\overline{V},\overline{E})$ 
from $(V,E)$ such that there exists a homeomorphism 
$\overline{H}:X_{(\overline{V},\overline{E})}\to X_{(V,E)}$ 
implementing an isomorphism between $AF(\overline{V},\overline{E})$ and 
the AF-subequivalence relation $\overline{R}$ of $R$ that we want. 
The Bratteli diagram $(\overline{V},\overline{E})$ will lend itself 
to apply Theorem \ref{glue} by, loosely speaking, 
transforming the subdiagrams $(W',F')_n$ in $(V,E)$ to replicas of $(W,F)$, 
thus making it possible 
to ``absorb'' the compact \'etale equivalence relation $K$. 
To construct $(\overline{V},\overline{E})$ from $(V,E)$ 
we first replace each vertex $v$ in $V_n$ (for $n\geq2$) 
which belong to the union of the vertex sets of 
$(W',F')_1,(W',F')_2,\dots,(W',F')_{n-1}$, 
by the vertices in $q_W^{-1}(\{v\})$. 
We retain the other vertices in $V_n$, and thus we get $\overline{V}_n$. 
We set $\overline{V}_0=V_0$ and $\overline{V}_1=V_1$. 
There is an obvious map $q_{\overline{V}}:\overline{V}\to V$, 
respecting gradings, which is surjective, 
and which can be considered to be 
an extension of the map $q_W:W\to W'$. 
Now we transfer in an obvious sense 
the subdiagram $(\widetilde{W},\widetilde{F})$ of $(V,E)$, 
again denoting it by $(\widetilde{W},\widetilde{F})$, 
and so $\overline{E}$ will contain the edge set $\widetilde{F}$. 
Similarly, in an obvious sense, 
we replace $(W',F')_n$ by $(W,F)_n$, $n\geq1$, 
where $(W,F)_n$ is defined in a similar way as we defined $(W',F')_n$. 
The new edge set $\overline{E}$ contains 
the collection of edges in $(W,F)_n$, $n\geq1$, and so, 
in particular, the edges $e_1,e_2,\dots$ lie in $\overline{E}$. 
(We will denote the path $(e_1,e_2,\dots)$ again by $x_\infty$.) 
Furthermore, 
if $e\in E$ is such that the vertices $i(e)$ and $t(e)$ do not lie in $L'$, 
then we retain $e$, and so $e\in\overline{E}$. 
Let $e\in E\setminus G'$, and let $v=i(e)$ and $w=t(e)$. 
If $v\in L'$, 
we replace $e$ by $\#(q_{\overline{V}}^{-1}(\{v\}))$ edges 
sourcing at each of the vertices $q_{\overline{V}}^{-1}(\{v\})$ 
and ranging at the same vertex, 
where this vertex can be chosen 
to be any vertex in $q_{\overline{V}}^{-1}(\{w\})$. 
If $v\notin L'$ and $w\in L'$, then 
we replace $e$ by an edge sourcing at $v$ 
and ranging at an arbitrary vertex in $q_{\overline{V}}^{-1}(\{w\})$. 
However, we require that 
the collection of these new edges will range 
at every vertex in $q_{\overline{V}}^{-1}(\{w\})$, 
as we consider all edges $e\in E\setminus G'$ 
such that $i(e)=v$ and $t(e)=w$, 
$v\in V_{n-1}$ and $w\in V_n$ being fixed ($n\geq1$). 
By (2) there are sufficiently many edges $e$ 
so that this can be achieved. 
So we have defined the edge set $\overline{E}$, 
and we have done this in such a way that 
$(\overline{V},\overline{E})$ is a simple Bratteli diagram 
by ensuring that for all $v\in\overline{V}_{n-1}$, $w\in\overline{V}_n$, 
there exists $e\in\overline{E}$ such that $i(e)=v$, $t(e)=w$. 
There is an obvious map $q_{\overline{E}}:\overline{E}\to E$, 
respecting gradings, which is surjective, 
and which can be considered to be an extension of the map $q_F:F\to F'$. 
By the definition of $q_{\overline{E}}$ and $q_{\overline{V}}$ 
we have that $i(q_{\overline{E}}(e))=q_{\overline{V}}(i(e))$, 
$t(q_{\overline{E}}(e))=q_{\overline{V}}(t(e))$ for $e\in\overline{E}$. 
Also, using Proposition \ref{transBratteli} (ii), 
it is easy to show that 
$q_{\overline{E}}:i^{-1}(\{v\})\to i^{-1}(\{q_{\overline{V}}(v)\})$ 
is a bijection for each $v\in\overline{V}$. 
We deduce that the map 
$(f_1,f_2,\dots,f_n)\mapsto
(q_{\overline{E}}(f_1),q_{\overline{E}}(f_2),\dots,q_{\overline{E}}(f_n))$ 
establishes a bijection between the finite paths of length $n$ 
from the top vertices of $(\overline{V},\overline{E})$ and $(V,E)$, 
respectively. 
We define the map 
$\overline{H}:X_{(\overline{V},\overline{E})}\to X_{(V,E)}$ by 
\[ x=(f_1,f_2,f_3,\dots)\mapsto 
\overline{H}(x)=(q_{\overline{E}}(f_1),q_{\overline{E}}(f_2),
q_{\overline{E}}(f_3),\dots). \]
It is now routinely checked that $\overline{H}$ is a homeomorphism, and 
that $\overline{H}\times\overline{H}$ maps $AF(\overline{V},\overline{E})$ 
isomorphically onto an AF-subequivalence relation 
$\overline{R}$ on $X_{(V,E)}$. 
Observe that 
$\overline{H}$ maps $X_{(W,F)_n}$ onto $X_{(W',F')_n}$ for all $n\geq1$. 
(In fact, by obvious identifications, 
this map is the same as $H:X_{(W,F)}\to X_{(W',F')}$ introduced earlier.) 
Also, $\overline{H}$ maps 
$X_{(\widetilde{W},\widetilde{F})}(\subset X_{(\overline{V},\overline{E})})$ 
onto $X_{(\widetilde{W},\widetilde{F})}(\subset X_{(V,E)})$, 
and $\overline{H}(x_\infty)=x_\infty$. 

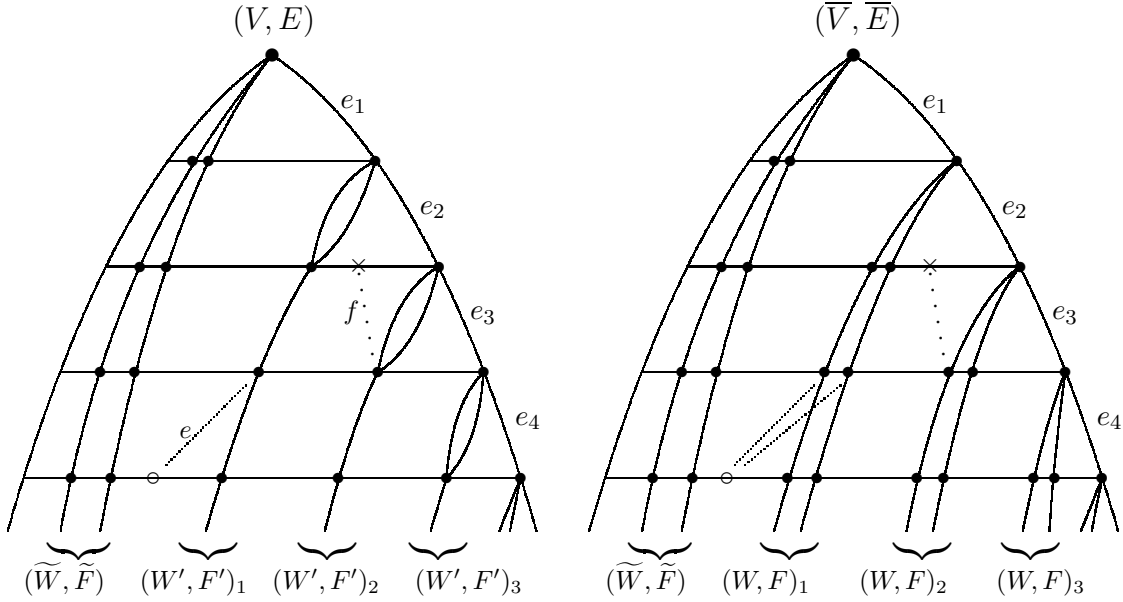
\begin{figure}
\begin{center}
\begin{picture}(420,240)

\put(100,210){\circle*{5}}
\qbezier(0,30)(40,170)(100,210)
\qbezier(200,30)(160,170)(100,210)
\put(85,220){$(V,E)$}

\put(61,170){\line(1,0){78}}
\put(37,130){\line(1,0){126}}
\put(20,90){\line(1,0){160}}
\put(6,50){\line(1,0){188}}

\put(139,170){\circle*{4}}
\put(163,130){\circle*{4}}
\put(180,90){\circle*{4}}
\put(194,50){\circle*{4}}

\put(126,190){\small{$e_1$}}
\put(156,150){\small{$e_2$}}
\put(175,110){\small{$e_3$}}
\put(192,70){\small{$e_4$}}

\qbezier(100,210)(40,140)(20,30)
\qbezier(100,210)(60,160)(35,30)
\put(15,26){$\underbrace{\qquad}$}
\put(6,8){\footnotesize{$(\widetilde{W},\widetilde{F})$}}
\put(70,170){\circle*{4}}\put(76,170){\circle*{4}}
\put(50,130){\circle*{4}}\put(60,130){\circle*{4}}
\put(35,90){\circle*{4}}\put(48,90){\circle*{4}}
\put(24,50){\circle*{4}}\put(39,50){\circle*{4}}

\qbezier(139,170)(120,160)(115,130)
\qbezier(139,170)(130,140)(115,130)
\qbezier(115,130)(95,100)(75,30)
\put(65,26){$\underbrace{\quad}$}
\put(50,8){\footnotesize{$(W',F')_1$}}
\put(115,130){\circle*{4}}
\put(95,90){\circle*{4}}
\put(81,50){\circle*{4}}

\qbezier(163,130)(146,120)(140,90)
\qbezier(163,130)(156,100)(140,90)
\qbezier(140,90)(130,70)(120,30)
\put(110,26){$\underbrace{\quad}$}
\put(100,8){\footnotesize{$(W',F')_2$}}
\put(140,90){\circle*{4}}
\put(125,50){\circle*{4}}

\qbezier(180,90)(168,78)(166,50)
\qbezier(180,90)(178,62)(166,50)
\qbezier(166,50)(163,40)(160,30)
\put(152,26){$\underbrace{\quad}$}
\put(154,8){\footnotesize{$(W',F')_3$}}
\put(166,50){\circle*{4}}

\qbezier(194,50)(190,40)(186,30)
\qbezier(194,50)(192,40)(190,30)

\put(320,210){\circle*{5}}
\qbezier(220,30)(260,170)(320,210)
\qbezier(420,30)(380,170)(320,210)
\put(305,220){$(\overline{V},\overline{E})$}

\put(281,170){\line(1,0){78}}
\put(257,130){\line(1,0){126}}
\put(240,90){\line(1,0){160}}
\put(226,50){\line(1,0){188}}

\put(359,170){\circle*{4}}
\put(383,130){\circle*{4}}
\put(400,90){\circle*{4}}
\put(414,50){\circle*{4}}

\put(346,190){\small{$e_1$}}
\put(376,150){\small{$e_2$}}
\put(395,110){\small{$e_3$}}
\put(412,70){\small{$e_4$}}

\qbezier(320,210)(260,140)(240,30)
\qbezier(320,210)(280,160)(255,30)
\put(235,26){$\underbrace{\qquad}$}
\put(226,8){\footnotesize{$(\widetilde{W},\widetilde{F})$}}
\put(290,170){\circle*{4}}\put(296,170){\circle*{4}}
\put(270,130){\circle*{4}}\put(280,130){\circle*{4}}
\put(255,90){\circle*{4}}\put(268,90){\circle*{4}}
\put(244,50){\circle*{4}}\put(259,50){\circle*{4}}

\qbezier(359,170)(320,140)(290,30)
\qbezier(359,170)(330,140)(300,30)
\put(285,26){$\underbrace{\quad}$}
\put(270,8){\footnotesize{$(W,F)_1$}}
\put(327,130){\circle*{4}}\put(334,130){\circle*{4}}
\put(309,90){\circle*{4}}\put(318,90){\circle*{4}}
\put(295,50){\circle*{4}}\put(306,50){\circle*{4}}

\qbezier(383,130)(356,110)(340,30)
\qbezier(383,130)(366,110)(350,30)
\put(335,26){$\underbrace{\quad}$}
\put(322,8){\footnotesize{$(W,F)_2$}}
\put(356,90){\circle*{4}}\put(365,90){\circle*{4}}
\put(344,50){\circle*{4}}\put(354,50){\circle*{4}}

\qbezier(400,90)(392,70)(384,30)
\qbezier(400,90)(396,60)(394,30)
\put(376,26){$\underbrace{\quad}$}
\put(374,8){\footnotesize{$(W,F)_3$}}
\put(388,50){\circle*{4}}\put(396,50){\circle*{4}}

\qbezier(414,50)(410,40)(406,30)
\qbezier(414,50)(412,40)(410,30)

\qbezier[30](90,85)(75,70)(60,55)
\put(55,50){\circle{4}}
\put(65,65){\footnotesize{$e$}}
\qbezier[30](305,85)(290,70)(275,55)
\qbezier[30](315,85)(297,70)(279,55)
\put(272,50){\circle{4}}

\put(129,128){\footnotesize{$\times$}}
\multiput(133,127)(1,-6){6}{\circle*{1}}
\put(127,110){\footnotesize{$f$}}
\put(345,128){\footnotesize{$\times$}}
\multiput(349,127)(1,-6){6}{\circle*{1}}

\end{picture}
\caption{Construction of the AF-subequivalence relation 
$\overline{R}\cong AF(\overline{V},\overline{E})$ of $R=AF(V,E)$, 
cf. Remark \ref{simplecase}}
\end{center}
\end{figure}

We claim that 
\[ R=\overline{R}\vee K_1\vee K_2\vee\dots, \tag{$*$} \]
where $K_n$, $n\geq1$, is the CEER on $X_{(W',F')_n}$ 
that corresponds to $K$ via the obvious map 
between $X_{(W',F')}$ and $X_{(W',F')_n}$. 
(We note for later use that, similarly, 
we have an obvious map between $X_{(W,F)}$ and $X_{(W,F)_n}$.) 
Clearly the right hand side of ($*$) is contained in the left hand side. 
Now $\overline{H}|_{(W,F)_n}$ implements 
an embedding of $AF((W,F)_n)(\cong AF(W,F)\cong R|_Y)$ 
into $AF((W',F')_n)(\cong AF(W',F')\cong(R|_Y)\vee K)$ 
which is transverse to $K_n$, and such that 
$AF((W',F')_n)\cong 
K_n\vee (\overline{H}\times\overline{H})(AF((W,F)_n))$ 
(cf. Proposition \ref{transBratteli}). 
To prove that ($*$) holds, 
let $(x,y)\in R$, 
and let $\overline{x},\overline{y}\in X_{(\overline{V},\overline{E})}$ 
such that $\overline{H}(\overline{x})=x,$ $\overline{H}(\overline{y})=y$. 
If the paths $x$ and $y$ agree from level $n$ on in $X_{(V,E)}$, 
then by the definition of $\overline{H}$ we must have that 
$q_{\overline{E}}(f_m)=q_{\overline{E}}(f'_m)$ and 
$q_{\overline{V}}(i(f_m))=q_{\overline{V}}(i(f'_m))$ for all $m>n$, 
where $\overline{x}=(f_1,f_2,\dots)$, $\overline{y}=(f'_1,f'_2,\dots)$. 
We first observe that 
$\overline{x}$ and $\overline{y}$ are cofinal paths 
in $X_{(\overline{V},\overline{E})}$, that is, 
$(\overline{x},\overline{y})\in AF(\overline{V},\overline{E})$, 
if and only if $i(f_m)=i(f'_m)$ for some $m>n$. 
Now it follows directly from the way we constructed 
the Bratteli diagram $(\overline{V},\overline{E})$ from $(V,E)$ 
that if $q_{\overline{E}}(f_m)(=q_{\overline{E}}(f'_m))$ does not 
belong to the edge set of a fixed $(W',F')_N$ for all $m>n$, 
then this situation will occur, and so 
$(\overline{H}(\overline{x}),\overline{H}(\overline{y}))
=(x,y)\in \overline{R}$. 
If on the other hand $q_{\overline{E}}(f_m)$ belongs to 
the edge set of a fixed $(W',F')_N$ for all $m>n$, 
then both $f_m$ and $f'_m$ belong to the edge set of 
a fixed $(W,F)_N$ for all $m>n$. 
We can then find paths $\tilde{x}$ and $\tilde{y}$ 
in $X_{(W,F)_N}(\subset X_{(\overline{V},\overline{E})})$ 
which are cofinal with $\overline{x}$ and $\overline{y}$, respectively, 
from level $n$ on. 
Hence we will have that 
$(\overline{H}(\tilde{x}),\overline{H}(\tilde{y}))\in AF((W',F')_N)$. 
Now $(\overline{x},\tilde{x}),(\overline{y},\tilde{y})
\in AF(\overline{V},\overline{E})$, and so 
$(\overline{H}(\overline{x}),\overline{H}(\tilde{x})),
(\overline{H}(\overline{y}),\overline{H}(\tilde{y}))\in \overline{R}$, 
i.e. 
$(x,\overline{H}(\tilde{x})),(y,\overline{H}(\tilde{y}))\in \overline{R}$. 
Combining all this we get that 
$(x,y)\in \overline{R}\vee K_N
\subset \overline{R}\vee K_1\vee K_2\vee\dots$, 
and so we have proved that ($*$) holds. 
(We remark that $\overline{R}$ is open in $R$, 
since it is a general fact that 
if $\overline{S}\subset S$ are \'etale equivalence relations on $X$, 
$\overline{S}$ having the relative topology from $S$, 
then $\overline{S}$ is open in $S$.) 

Let $(\overline{L},\overline{G})$ be the subdiagram of 
$(\overline{V},\overline{E})$ whose edge set consists of $\widetilde{F}$ 
and the union of the edge sets belonging to 
$(W,F)_1,(W,F)_2,\dots$, 
and let $(\overline{L'},\overline{G'})$ be the subdiagram of 
$(\overline{V},\overline{E})$ whose edge set consists of 
the union of the edge sets belonging to $(W,F)_1,(W,F)_2,\dots$. 
(Note that $(\overline{L},\overline{G})$ and $(\overline{L'},\overline{G'})$ 
are analogous to the subdiagrams $(L,G)$ and $(L',G')$, 
respectively, of $(V,E)$ that we defined above.) 
By (2) it follows that 
$(\overline{L'},\overline{G'})$ is a thin subdiagram. 
Also, it is easily seen that 
$(\widetilde{W},\widetilde{F})$ is thin 
in $(\overline{V},\overline{E})$, 
and hence $(\overline{L},\overline{G})$ is thin 
in $(\overline{V},\overline{E})$. 
(Observe that 
both $X_{(\overline{L},\overline{G})}$ and $X_{(\overline{L'},\overline{G'})}$ 
are homeomorphic to $(Y\times\N)\cup\{x_\infty\}$, 
where $\N=\{1,2,3,\dots\}$ (with discrete topology), 
and $x_\infty$ is the point at infinity 
of the one-point compactification of $Y\times\N$.) 
Clearly $Z=X_{(\overline{L},\overline{G})}$ and 
$Z'=X_{(\overline{L'},\overline{G'})}$ are 
$\widetilde{R}$-\'etale closed subsets of $X_{(\overline{V},\overline{E})}$ 
(where we for convenience write $\widetilde{R}$ 
for $AF(\overline{V},\overline{E})$). 
We will define a homeomorphism $\alpha:Z\to Z'$ 
which implements an isomorphism 
between $\widetilde{R}|_Z(\cong AF(\overline{L},\overline{G}))$ 
and $\widetilde{R}|_{Z'}(\cong AF(\overline{L'},\overline{G'}))$. 
At the same time 
we want $\alpha\times\alpha$ to map $K$ 
(on $Y=X_{(\widetilde{W},\widetilde{F})}
\subset X_{(\overline{V},\overline{E})}$) isomorphically to 
$K_1$ (on $X_{(W,F)_1}$), 
and $K_n$ (on $X_{(W,F)_n}$) isomorphically to 
$K_{n+1}$ (on $X_{(W,F)_{n+1}}$) for all $n\geq1$. 
(Retaining the previous notation, 
we mean by $K_n$ on $X_{(W,F)_n}$ 
the CEER corresponding to $K_n$ on 
$X_{(W',F')_n}$---$\overline{H}^{-1}\times\overline{H}^{-1}$ 
transferring $K_n$ on $X_{(W',F')_n}$ to $K_n$ on $X_{(W,F)_n}$.) 
There is an obvious way 
to define $\alpha:X_{(W,F)_n}\to X_{(W,F)_{n+1}}$ 
satisfying our requirements, 
using the fact that 
$(W,F)_n$ and $(W,F)_{n+1}$ are essentially replicas of $(W,F)$. 
Now consider $Y=X_{(\widetilde{W},\widetilde{F})}$, 
where $X_{(\widetilde{W},\widetilde{F})}
\subset X_{(\overline{V},\overline{E})}$. 
We have $AF(\widetilde{W},\widetilde{F})\cong\widetilde{R}|_Y$ 
(which clearly may be identified with $R|_Y$), and 
$\widetilde{R}|_Y\cong\widetilde{R}|_{X_{(W,F)_1}}(\cong AF(W,F))$. 
Furthermore, 
$(R|_Y)\vee K\cong(\widetilde{R}|_{X_{(W,F)_1}})\vee K_1$. 
The various isomorphism maps are naturally related, 
and it follows that 
there exists $\alpha:X_{(\widetilde{W},\widetilde{F})}\to X_{(W,F)_1}$ 
satisfying our requirements. 
(We omit the easily checked details.) 
Furthermore, note that 
$\alpha$ implements an isomorphism 
between $(\widetilde{R}|_Y)\vee K$ and 
$(\widetilde{R}|_{X_{(W,F)_1}})\vee K_1$. 
Also, $\overline{H}\times\overline{H}$ implements an isomorphism 
between the latter and $R|_{X_{(W',F')_1}}$. 
We define eventually $\alpha:Z\to Z'$ 
by patching together the various $\alpha$'s above, 
letting $\alpha(x_\infty)=x_\infty$, 
and it is straightforward to verify that 
this $\alpha$ satisfies all our requirements. 
By Theorem \ref{glue} 
there exists an extension 
$\overline{\alpha}:X_{(\overline{V},\overline{E})}
\to X_{(\overline{V},\overline{E})}$ of $\alpha$ 
which implements an automorphism of 
$\widetilde{R}=AF(\overline{V},\overline{E})$. 
Let $h=\overline{H}\circ\overline{\alpha}\circ\overline{H}^{-1}$. 
Then $h$ is a homeomorphism on $X=X_{(V,E)}$. 
Recall that 
$\overline{H}\times\overline{H}(\widetilde{R})=\overline{R}$. 
By ($*$) we have that 
\[ R\vee K=\overline{R}\vee K\vee K_1\vee K_2\vee\dots. \tag{$**$} \]
We note that 
\begin{align*}
h\times h(\overline{R})
&=(\overline{H}\times\overline{H})
\circ(\overline{\alpha}\times\overline{\alpha})
\circ(\overline{H}^{-1}\times\overline{H}^{-1})(\overline{R}) \\
&=(\overline{H}\times\overline{H})
\circ(\overline{\alpha}\times\overline{\alpha})(\widetilde{R}) \\
&=(\overline{H}\times\overline{H})(\widetilde{R}) \\
&=\overline{R}. 
\end{align*}
Also, 
\begin{align*}
h\times h(K)
&=(\overline{H}\times\overline{H})
\circ(\overline{\alpha}\times\overline{\alpha})
\circ(\overline{H}^{-1}\times\overline{H}^{-1})(K) \\
&=(\overline{H}\times\overline{H})
\circ(\overline{\alpha}\times\overline{\alpha})(K) \\
&=(\overline{H}\times\overline{H})(K_1) \\
&=K_1, 
\end{align*}
where we make the obvious identifications, referred to above, 
with $K$ (respectively $K_1$) 
as CEERs on $X_{(\widetilde{W},\widetilde{F})}(\subset X_{(V,E)})$ 
and $X_{(\widetilde{W},\widetilde{F})}
(\subset X_{(\overline{V},\overline{E})})$ 
(respectively, $X_{(W',F')_1}(\subset X_{(V,E)})$ 
and $X_{(W,F)_1}(\subset X_{(\overline{V},\overline{E})})$). 
Similarly we show that 
$h\times h(K_n)=K_{n+1}$ for $n\geq1$. 
By ($**$) we get 
\begin{align*}
h\times h(R\vee K)
&=(h\times h)(\overline{R})\vee(h\times h)(K)\vee(h\times h)(K_1)\vee\dots \\
&=\overline{R}\vee K_1\vee K_2\vee\dots \\
&=R. 
\end{align*}
This completes the proof of the main assertion, (i), of the theorem. 

Assertions (ii) and (iii) are now immediately clear 
since $h(Y)=X_{(W',F')_1}$ and 
$h|_Y$ (by definition of $\overline{\alpha}|_Y$) implements 
an isomorphism between $(R|_Y)\vee K$ and $R|_{h(Y)}(\cong AF((W',F')_1))$. 
This finishes the proof of the theorem. 
\end{proof}

\begin{rem}\label{simplecase}
In order to understand 
the idea behind the proof of the absorption theorem better, 
it is instructive to look at the simplest (non-trivial) case, 
namely when $Y=\{y_1,y_2\}$ consists of two points $y_1$ and $y_2$, 
such that $(y_1,y_2)\notin R$. 
(Even in this simple case, 
the conclusion one can draw from the absorption theorem 
is highly non-trivial.) 
The compact \'etale equivalence relation $K$ on $Y$ 
that is transverse to $R|_Y(=\Delta_Y)$ is the following: 
$K=\Delta_Y\cup\{(y_1,y_2),(y_2,y_1)\}$. 
The Bratteli diagrams $(\widetilde{W},\widetilde{F})$ and $(W,F)$ 
for $R|_Y$ are trees consisting of two paths with no vertices in common, 
except the top one. 
The Bratteli diagram $(W',F')$ for $(R|_Y)\vee K$ starts 
with two edges forming a loop, 
and then a single path. 
In Figure 6 the scenario in this case is illustrated. 
(We have indicated 
how one constructs the new Bratteli diagram 
$(\overline{V},\overline{E})$ from $(V,E)$ (where $R=AF(V,E)$) 
by exhibiting two specific edges $e$ and $f$ in $E$ 
and how they give rise to new edges 
$q^{-1}_{\overline{E}}(\{e\})$ and $q^{-1}_{\overline{E}}(\{f\})$, 
respectively, in $\overline{E}$.) 
The example considered here corresponds to 
a very special case of transversality 
arising from an action of the (finite) group $G=\Z/2\Z=\{0,1\}$, 
where $\alpha_1:Y\to Y$ sends $y_1$ to $y_2$ and $y_2$ to $y_1$ 
(cf. the comments after Definition \ref{transverse}). 
In the paper \cite{GPS2} 
the general $\Z/2\Z$-action case is considered 
(even though it is formulated slightly different there). 
\end{rem}

\begin{rem}\label{addition}
Theorem \ref{main} is the key technical ingredient 
in the proof of the following result: 
\end{rem}

\begin{theo}[\cite{GMPS}]
Every minimal action of $\Z^2$ on a Cantor set is 
orbit equivalent to an AF-equivalence relation, 
and consequently also orbit equivalent to a minimal $\Z$-action. 
\end{theo}


\begin{thebibliography}{ABCD}
\bibitem[GMPS]{GMPS}
T. Giordano, H, Matui, I. F. Putnam and C. F. Skau, 
\textit{Orbit equivalence for Cantor minimal $\Z^2$-systems}, 
preprint. math.DS/0609668. 
\bibitem[GPS1]{GPS1}
T. Giordano, I. F. Putnam and C. F. Skau, 
\textit{Topological orbit equivalence and $C^*$-crossed products}, 
J. Reine  Angew. Math. 469 (1995), 51--111. 
\bibitem[GPS2]{GPS2}
T. Giordano, I. F. Putnam and C. F. Skau, 
\textit{Affable equivalence relations and orbit structure 
of Cantor dynamical systems}, 
Ergodic Theory Dynam. Systems 24 (2004), 441--475. 
\bibitem[M]{M}
M. Molberg, 
\textit{AF-equivalence relations}, 
Math. Scand. 99 (2006), 247--256. 
\bibitem[P]{P}
A. L. T. Paterson, 
\textit{Groupoids, inverse semigroups, and their operator algebras}, 
Progress in Mathematics, 170. 
Birkh\"auser Boston, Inc., Boston, MA, 1999. 
\end{thebibliography}
\end{document}